\documentclass[review]{elsarticle}

\usepackage{lineno,hyperref}
\modulolinenumbers[5]
\journal{Journal of XXX}

\usepackage{amsmath,amssymb,amsthm,amsfonts,afterpage,float,bm,stmaryrd,paralist,wrapfig,bbm,soul}

\usepackage{cancel}
\usepackage{color}
\usepackage{fancyhdr}
\usepackage{graphics}
\usepackage{hyperref}
\usepackage{graphicx,epsf}
\usepackage{makeidx}
\usepackage{epsfig}
\usepackage{lscape}

\usepackage{empheq,xcolor}

\pssilent

\newtheorem{theorem}{Theorem}[section]
\newtheorem{lemma}{Lemma}[section]
\newtheorem{proposition}{Proposition}[section]
\newtheorem{corollary}{Corollary}[section]
\newtheorem{conjecture}{Conjecture}[section]

\newtheorem{remark}{Remark}[section]
\newtheorem{example}{Example}[section]

\numberwithin{equation}{section}
\numberwithin{figure}{section}
\numberwithin{table}{section}

\def\XXint#1#2#3{{\setbox0=\hbox{$#1{#2#3}{\int}$}
\vcenter{\hbox{$#2#3$}}\kern-.51\wd0}}

\renewcommand{\div}{\operatorname{div}}

\graphicspath{{Figure/}}

\begin{document}

\setlength{\pdfpageheight}{\paperheight}
\setlength{\pdfpagewidth}{\paperwidth}
\title{Nonlocal Effect on a Generalized Ohta-Kawasaki Model}
\author{Wangbo Luo}
\author{Yanxiang Zhao}
\address{Department of Mathematics, George Washington University, Washington D.C., 20052}

\begin{abstract}
We propose a nonlocal Ohta-Kawasaki model to study the nonlocal effect on the pattern formation of some binary systems with general long-range interactions. While the nonlocal Ohta-Kawasaki model displays similar bubble patterns as the standard Ohta-Kawasaki model, by performing Fourier analysis, we find that the optimal number of bubbles for the nonlocal model may have an upper bound no matter how large the repulsive strength is. The existence of such an upper bound is characterized by the eigenvalues of the nonlocal kernels. Additionally we explore the conditions under which the nonlocal horizon parameter may promote or demote the bubble splitting, and apply the analysis framework to several case studies for various nonlocal operators.
\end{abstract}


\date{\today}
\maketitle


\section{Introduction}\label{sec:Introduction}

Ohta-Kawasaki (OK) model is introduced in \cite{OhtaKawasaki_Macromolecules1986} and has been extensively applied for the study of phase separation of diblock copolymers, which
have generated much interest in materials science in the past years due to their remarkable ability for self-assembly into nanoscale ordered structures \cite{Hamley_Wiley2004}. Diblock copolymers are chain molecules made by two different segment species, say $A$ and $B$. Due to the chemical incompatibility, the two species tend to be phase-separated; on the other hand, the two species are connected by covalent chemical bonds, which leads to the so-called microphase separation. 

Block copolymers provide simple and easily controlled materials for the study of self-assembly. Mean field theories, with an associated free energy functional, have proven very useful in the understanding and prediction of the pattern morphology \cite{Bates_PhysToday1999, Choksi_NonlinearScience2001}; Experimentally and numerically, it is well known that block copolymer system displays periodic structures such as lamellar, spherical, bicontinuous gyroids \cite{BahianaOono_PRA1990, Bates_PhysToday1999,HasegawaTannakaYamasakiHashimoto_Macromolecules1987, ZhengWang_Micromolecules1995,TangQiuYang_PRE2004,LiJiangChen_SoftMatter2013}. However some other interesting microphases are overlooked and have not been well studied theoretically. Unusual patterns were found, such as spherical/square bubbles mixing in square lattice, elliptical bubbles in hexagonal lattice \cite{HeZhouKanLiang_JCP2015} and elongated hexagons \cite{LiuLiQiu_Macromolecules2012, LeeRyuetal_ACSNano2018}. These asymmetric microphases cannot be predicted by the usual symmetric Ohta-Kawasaki/Ohta-Nakazawa theory. Besides, quantitative study as a comparison to theoretical study \cite{Choksi_NonlinearScience2001} has not been well addressed. Therefore one needs to examine the variational problem in its full generality from a mathematically more sophisticated point of view, one which in particular allows for a fuller analysis of the competition between the terms in the energy. Therefore we propose a Nonlocal Ohta-Kawasaki (NOK) model by a free energy functional:
\begin{align}\label{functional:NOK}
E^{\text{NOK}}_{\epsilon}[u] = \int_{\mathbb{T}^d} \dfrac{\epsilon}{2}|\nabla u|^2 + \dfrac{1}{\epsilon}W(u)\ \text{d}x + \dfrac{\gamma}{2}\int_{\mathbb{T}^d} |(\mathcal{L}_{\delta})^{-\frac{1}{2}}(u-\omega)|^2\ \text{d}x,
\end{align}
with a possible volume constraint
\begin{align}\label{eqn:Volume}
\int_{\mathbb{T}^d} (u - \omega)\ \text{d}x = 0.
\end{align}
Here $\mathbb{T}^d = \prod_{i=1}^d [-X_i, X_i] \subset \mathbb{R}^d, d =1, 2, 3$ denotes a periodic box and $0<\epsilon \ll 1$ is an interface parameter that indicates the system is in deep segregation regime. $u = u(x)$ is a phase field labeling function which represents the concentration of species $A$. By the assumption of incompressibility for the binary system, the concentration of species $B$ can be implicitly represented by $1-u(x)$. Function $W(u) = 18(u^2-u)^2$ is a double well potential which enforces the phase field function $\phi$ to be equal to 1 inside the interface and 0 outside the interface. Near the interfacial region, the phase field function $\phi$ rapidly but smoothly transitions from 0 to 1. The first integral in (\ref{functional:NOK}) is a local surface energy which represents the short-range interaction between the chain molecules and favors the large domain; while the second integral in (\ref{functional:NOK}) is a term for the long-range (nonlocal) repulsive interaction with $\gamma >0$ being the strength of the repulsive force. Finally, $\omega\in(0,1)$ is the relative volume of the species $A$. Since species $A$ and $B$ are incompressible, it is sufficient to consider $\omega \in (0,1/2)$ as otherwise one can simply swap the roles of $A$ and $B$.

$\mathcal{L}_{\delta}$ is a positive semi-definite/definite operator $\mathcal{L}_{\delta}: L^2_{\text{per}}(\mathbbm{T}^d) \rightarrow L^2_{\text{per}}(\mathbbm{T}^d)$, defined as
\begin{align}\label{kernel_Lij}
\textstyle{
    \mathcal{L}_{\delta}: u(x) \rightarrow \int_{\mathbb{T}^d} L_{\delta}(x-y)(u(x)-u(y))\ \text{d}y,
}    
\end{align}
where $L_{\delta}$ is a kernel function and assumed to be nonnegative, radially symmetric ($L_{\delta}(s) = L_{\delta}(|s|)$ for any $s$), compactly supported and has a bounded second moment \cite{Du_Book2020}. The horizon $\delta>0$ is introduced to measure the range of nonlocal interaction by $\mathcal{L}_{\delta}$. The value of $\delta$ is usually small and is restricted in $(0,\pi]$ in this work.  For a comprehensive review on nonlocal modeling, analysis and computation, we refer the interested readers to a recent monograph \cite{Du_Book2020} by Du.

 The previous work is mostly restricted to radially symmetric kernels \cite{Du_Book2020}. In the model (\ref{functional:NOK}), a more general kernel $L_{\delta}(s)$ can be incorporated for the consideration. For instance, $L_{\delta}$ can be with compactly supported in a square domain and satisfy $L_{\delta}(\|s\|_1) = L_{\delta}(s)$ or $L_{\delta}(\|s\|_{\infty}) = L_{\delta}(s)$. Note that the ``nonlocality" in $E_{\epsilon}^{\text{NOK}}[u]$ is twofold: one is due to the inverse $(\mathcal{L}_{\delta})^{-1}$ of the operator $\mathcal{L}_{\delta}$; the other is due to the nature of nonlocality of $\mathcal{L}_{\delta}$ itself, characterized by the nonlocal kernel $L_{\delta}$. By considering a general kernel $L_{\delta}$ (for instance, non-radial kernels), one can investigate various types of global minimizers for NOK model, some of which may display different patterns from those generated by the original OK model. 

Since it is proposed by Ohta and Kawasaki in \cite{OhtaKawasaki_Macromolecules1986}, there has been extensive work on theoretical analysis and numerical methods for the original Ohka-Kawasaki(OK) model, which is the free energy functional (\ref{functional:NOK}) by replacing $\mathcal{L}_{\delta}$ by $-\Delta$. For instance, in \cite{RenWei_SIAM2000, RenTruskinovsky_Elasticity2000}, the authors characterized the minimizers of OK model, and considered some other nonlocal geometric variational problem involving $(I-\gamma^2\Delta)^{-1}$, which can lead to the Fitzhugh-Nagumo system \cite{ChenChoiHuRen_SIMA2018}. Choski in \cite{Choksi_QAM2012} performed asymptotic analysis for the global minimizers of OK model. Recently Chan, Nejad and Wei \cite{ChanNejadWei_PhysicaD2019} considered a variant of OK model in which $\nabla u$ is replaced by a fractional diffusion. In their work, the authors prove the $\Gamma$-convergence and the existence of the global minimizers, and provide an energy growth estimate for their model. In recent years, some numerical schemes were also developed for OK model. For example, \cite{Benesova_SINA2014} studies an implicit midpoint spectral approximation for the equilibrium of OK model. \cite{ChengYangShen_JCP2017} adopts the IEQ method to study the diblock copolymer model. However, the existing works mainly focus on the Cahn-Hilliard dynamics, namely, the $H^{-1}$ gradient flow dynamics of OK model. Recently efforts have been made to design numerical schemes for the $L^2$ gradient flow dynamics of OK model, such as operator-splitting energy stable methods \cite{XuZhao_JSC2019} and maximum principle preserving methods \cite{XuZhao_JSC2020}.

The main contribution of this paper lies in several aspects. Firstly, a nonlocal OK model is proposed which includes a general long-range interaction term induced by the inverse of a nonlocal operator $\mathcal{L}_{\delta}$. The inclusion of $\mathcal{L}_{\delta}$ in the NOK model (\ref{functional:NOK}) can characterize a broader class of features of microphase separation and provide more insights on theoretical studies of these subjects. Secondly, when considering a special feasible set in the one-dimensional case, bubble functions of equal size and equal distance, we find that the optimal number of bubbles, as a function of the long-range repulsive force $\gamma$, may have an upper bound. This is in contrast to the result of OK model, in which the optimal number of bubbles grows to $\infty$ as $\gamma\rightarrow\infty$. Thirdly, we explore the $\delta$-effect on the optimal number of bubbles. Under some mild conditions, we can perform analysis for the effect of the nonlocal horizon $\delta$ on the optimal number of bubbles, namely, whether $\delta$ promotes or demotes the bubble splitting.

The rest of the paper is organized as follows. Section \ref{section:Sharp} discusses the sharp interface limit of the NOK model. In section \ref{section:Nonlocal}, we study the nonlocal effect on the minimizers of the NOK model. More specifically, we characterize the main features of the minimizers in subsections \ref{subsection:AOmega} and \ref{subsection: UOmega}. Then we study the nonlocal effect of $\gamma$ on the minimizers of NOK model for power kernels in subsection \ref{subsection:PowerKernel} and Gauss-type kernel in subsection \ref{subsection:GaussKernel}. In subsection \ref{subsection:delta}, $\delta$-effect on minimizers is considered under some mild conditions, followed by several case studies for various types of long-range interaction. Finally the conclusion is drawn and several directions of the future work are discussed in section \ref{section:Conclusion}.

\section{Sharp Interface Limit}\label{section:Sharp}

In this section, we present that NOK model (\ref{functional:NOK}) has a $\Gamma$-limit $E_0^{\text{NOK}}$. To this end, we define the NOK functional more rigorously as follows \cite{RenWei_SIAM2000}
\begin{align}\label{functional:NOK02}
E^{\text{NOK}}_{\epsilon}[u] =
\begin{cases}
 \int_{\mathbb{T}^d} \dfrac{\epsilon}{2}|\nabla u|^2 + \dfrac{1}{\epsilon}W(u)\ \text{d}x + \dfrac{\gamma}{2}\int_{\mathbb{T}^d} |(\mathcal{L}_{\delta})^{-\frac{1}{2}}(f(u)-\omega)|^2\ \text{d}x, \quad &u\in L_{\text{per},\omega}^2(\mathbbm{T}^d)\cap W^{1,2}(\mathbbm{T}^d); \\
 \infty, \quad &u\in L_{\text{per},\omega}^2(\mathbbm{T}^d)\backslash W^{1,2}(\mathbbm{T}^d).
 \end{cases}
\end{align}
Here the function space $L_{\text{per},\omega}^2(\mathbbm{T}^d)$ is the set of all squared integrable periodic functions satisfying the volume condition (\ref{eqn:Volume}). On the other hand, we define 
\begin{align}\label{functional:NOK03}
E^{\text{NOK}}_{0}[u] =
\begin{cases}
\|Du\|({\mathbbm{T}^d}) + \dfrac{\gamma}{2}\int_{\mathbb{T}^d} |(\mathcal{L}_{\delta})^{-\frac{1}{2}}(f(u)-\omega)|^2\ \text{d}x, \quad &u\in L_{\text{per},\omega}^2(\mathbbm{T}^d)\cap \text{BV}(\mathbbm{T}^d,\{0,1\}); \\
 \infty, \quad &u\in L_{\text{per},\omega}^2(\mathbbm{T}^d)\backslash \text{BV}(\mathbbm{T}^d,\{0,1\}).
 \end{cases}
\end{align}
Here $\text{BV}(\mathbbm{T}^d,\{0,1\}): = \{u\in\text{BV}(\mathbbm{T}^d): u(x) = 0 \text{\ or\ } 1 \text{\ for\ a.e.\ } x\in\mathbbm{T}^d\}$. $\text{BV}(\mathbbm{T}^d)$ is the space of functions of bounded variation. $\|Du\|$ is the absolute value of the distributional derivative of $u$, which can be viewed as a finite nonnegative measure on $\mathbbm{T}^d$. $\|Du\|({\mathbbm{T}^d})$ is the size of ${\mathbbm{T}^d}$ under this measure.

\begin{proposition}[$\Gamma$-convergence of $E^{\text{NOK}}_{\epsilon}$ to $E^{\text{NOK}}_{0}$]
For any $\{u_{\epsilon}\}\subset L_{\emph{per},\omega}^2(\mathbbm{T}^d)$ such that $\lim_{\epsilon\rightarrow0}u_{\epsilon} = u_0$, we have
\[
\liminf_{\epsilon\rightarrow0} E^{\emph{NOK}}_{\epsilon}(u_{\epsilon}) \ge E^{\emph{NOK}}_{0}(u_{0});
\]
for any $u\in L_{\emph{per},\omega}^2(\mathbbm{T}^d)\cap \emph{BV}(\mathbbm{T}^d,\{0,1\})$, there exists a family $\{u_{\epsilon}\}\subset L_{\emph{per},\omega}^2(\mathbbm{T}^d)$ such that $\lim_{\epsilon\rightarrow0}u_{\epsilon} = u_0$, and 
\[
\liminf_{\epsilon\rightarrow0} E^{\emph{NOK}}_{\epsilon}(u_{\epsilon}) \le E^{\emph{NOK}}_{0}(u_{0}).
\]
\end{proposition}

\begin{proof}
The proof can be completed by repeating the work in Modica \cite{Modica_ARMA1987} and Ren and Wei \cite{RenWei_SIAM2000}. The only condition we need to verify is the continuity of the functional
\begin{align}\label{eqn:continuity}
\dfrac{\gamma}{2}\int_{\mathbb{T}^d} |(\mathcal{L}_{\delta})^{-\frac{1}{2}}(u-\omega)|^2\ \text{d}x, 
\end{align}
or equivalently the boundedness of $\|\mathcal{L}_{\delta}^{-1}\|_{2}$. Indeed, one can apply the triangle inequality to have
\[
\|\mathcal{L}_{\delta}^{-1}\|_{2} \le \|\mathcal{L}_{\delta}^{-1} - \mathcal{L}_0^{-1}\|_{2} + \| \mathcal{L}_0^{-1}\|_{2} = \text{I} + \text{II}.
\]
The estimate I was proved in Lemma 1 of \cite{DuYang_SINA2016}, and estimate II is a standard result of elliptic regularity, see for instance \cite{Adams_SobolevSpace}. Then it follows that
\[
\|\mathcal{L}_{\delta}^{-1}\|_{2} \le  C_1\delta^2 + C_2
\]
in which $C_1$ and $C_2$ are two generic constants only dependent of the domain $\mathbb{T}^d$. Therefore the continuity of the functional (\ref{eqn:continuity}) is verified, and we have the following $\Gamma$-convergence hold
\[
E^{\text{NOK}}_{\epsilon} \overset{\Gamma}{\rightarrow} E^{\text{NOK}}_{0}.
\]
Here the notation of $\Gamma$-convergence is defined by the two inequalities in this proposition.
\end{proof}

When $u = \chi_{\Omega}$ for some Lebesgue measurable set $\Omega\subseteq \mathbbm{T}^d$, the energy $E_0^{\text{NOK}}$ of (\ref{functional:NOK03}) in the $\Gamma$-limit can be recast as 
\begin{align}\label{eqn:SharpLimit}
E_0^{\text{NOK}}[\chi_{\Omega}] = P_{\mathbb{T}^d}(\Omega) + \dfrac{\gamma}{2}\int_{\mathbb{T}^d} |(\mathcal{L}_{\delta})^{-\frac{1}{2}}( \chi_{\Omega} - \omega)|^2\ \text{d}x, 
\end{align}
where $P_{\mathbb{T}^d}(\Omega): = \int_{\mathbbm{T}^d}|\nabla\chi_{\Omega}|dx$, and 
\[
\int_{\mathbbm{T}^d} |\nabla f| dx : = \sup \left\{ \int_{\mathbbm{T}^d} f \div g \ dx: g \text{\ is\ compactly\ supported\ and\ } C^1, |g|\le 1 \text{\ in\ } \mathbbm{T}^d\right\}. 
\]
In other words, $P_{\mathbb{T}^d}(\Omega)$ stands for the perimeter of $\Omega$ in $\mathbbm{T}^d$. Hereafter, we will use $E_0^{\text{NOK}}$ in (\ref{eqn:SharpLimit}) as the sharp interface formulation of the NOK functional for the characteristic function of some Lebesgue measurable set $\Omega$.

\section{Nonlocal Effect on Minimizers of NOK Model}\label{section:Nonlocal}

From now on, we will focus on the one-dimensional NOK model. To provide the sharp interface formulation in the one-dimensional case, we introduce the following notations.

For a fixed positive integer $N$, let $\mathcal{A}_N^{\omega}$ be the set of periodic step functions of the form
\begin{align}\label{eqn:stepfunction}
u(x) = \sum_{i=1}^{2N} \frac{1+(-1)^{i+1}}{2} \chi_{[x_i,x_{i+1})}(x), x\in[0,2\pi)
\end{align}
with $0 = x_1 < x_2 < \cdots < x_{2N}<x_{2N+1} = 1$, such that 
\begin{align}\label{eqn:stepfunction_constraint}
\frac{1}{2\pi}\int_0^{2\pi} u(x) dx = \omega, \quad \text{i.e.}\quad (x_2 - x_1) + \cdots + (x_{2N}-x_{2N-1}) = 2\pi\omega.
\end{align}
Note that $x_1 = 0$ is specified owing to the periodicity of $u(x)$. Geometrically, $\mathcal{A}_N^{\omega}$ represents the set of periodic $\{0,1\}$-step functions with $N$ bubbles (where $u(x) = 1$) of relative total area $\omega$. Hereafter we call a function $u(x) \in \mathcal{A}_N^{\omega}$ as a $N$-bubble function. We will further define $U_N^{\omega}$ as the bubble function in $\mathcal{A}_N^{\omega}$ such that the $N$ bubbles are of equal size and equal distance, namely, $U_N^{\omega}$ is of the form (\ref{eqn:stepfunction}) with
\begin{align}
x_{2i-1} = \frac{2\pi i}{N}, x_{2i} = x_{2i-1} + \frac{2\pi \omega}{N}, \quad i=1,\cdots,N.
\end{align}
We will also define $\mathcal{U}^{\omega}$ as the set of bubble functions of equal size and equal distance,
\begin{align}\label{eqn:stepfunction_equalsize_equaldistance}
\mathcal{U}^{\omega} : = \{U_N^{\omega}\}_{N=1}^{\infty}.
\end{align}

With the notations introduced above, we can have the the sharp interface formulation (\ref{eqn:SharpLimit}) reduces to 1D case as
\begin{align}\label{eqn:SharpLimit_1d}
E_0^{\text{NOK}}[u(x)] = 2N + \dfrac{\gamma}{2}\int_0^{2\pi} |(\mathcal{L}_{\delta})^{-\frac{1}{2}}(u(x)-\omega)|^2\ \text{d}x,\quad u(x) \in \mathcal{A}_N^{\omega}.
\end{align}
Minimizing the 1D problem (\ref{eqn:SharpLimit_1d}) in the sharp interface limit can be treated as a two-step process \cite{RenTruskinovsky_Elasticity2000}. The first step is to minimize $E_0^{\text{NOK}}$ in $\mathcal{A}_N^{\omega}$. Once the minimizer is characterized, say $A_N^{\omega}$, one can move to the second step to minimize $E_0^{\text{NOK}}[A_N^{\omega}]$ over $N\in\mathbb{Z}^+$.

\begin{figure}[t] 
\begin{center}
\includegraphics[width=0.4\linewidth]{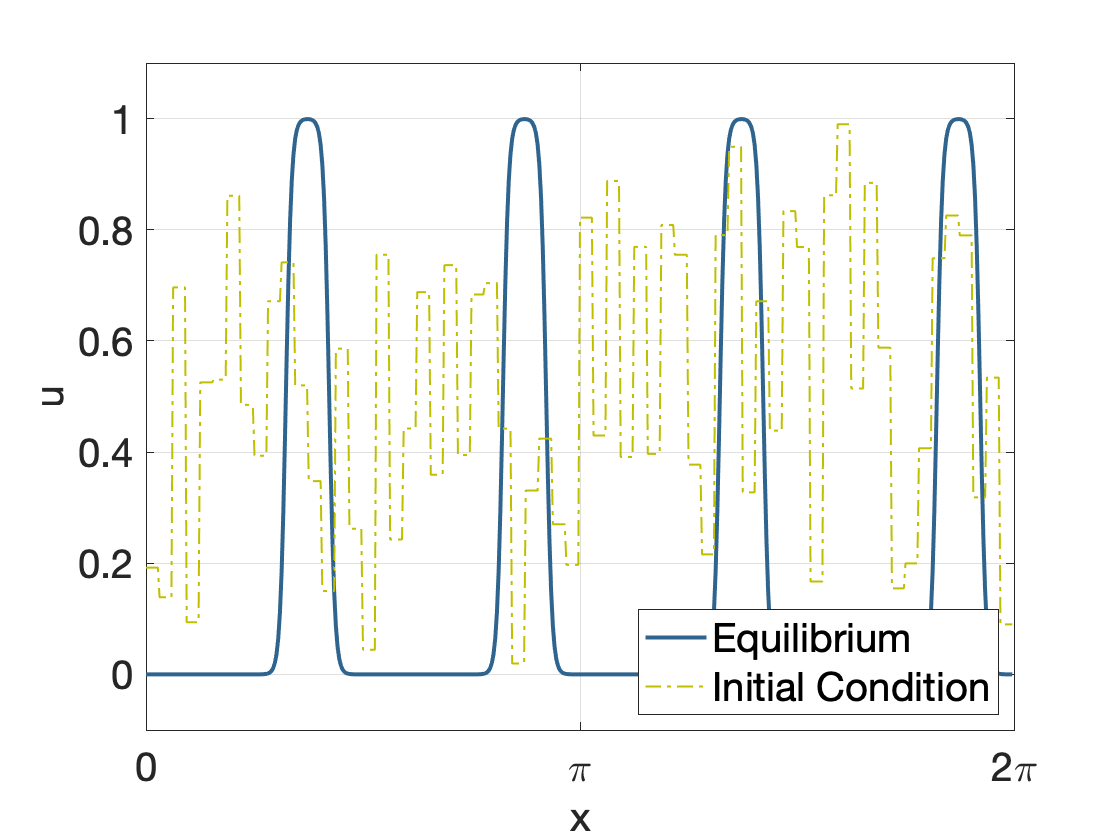}
\includegraphics[width=0.4\linewidth]{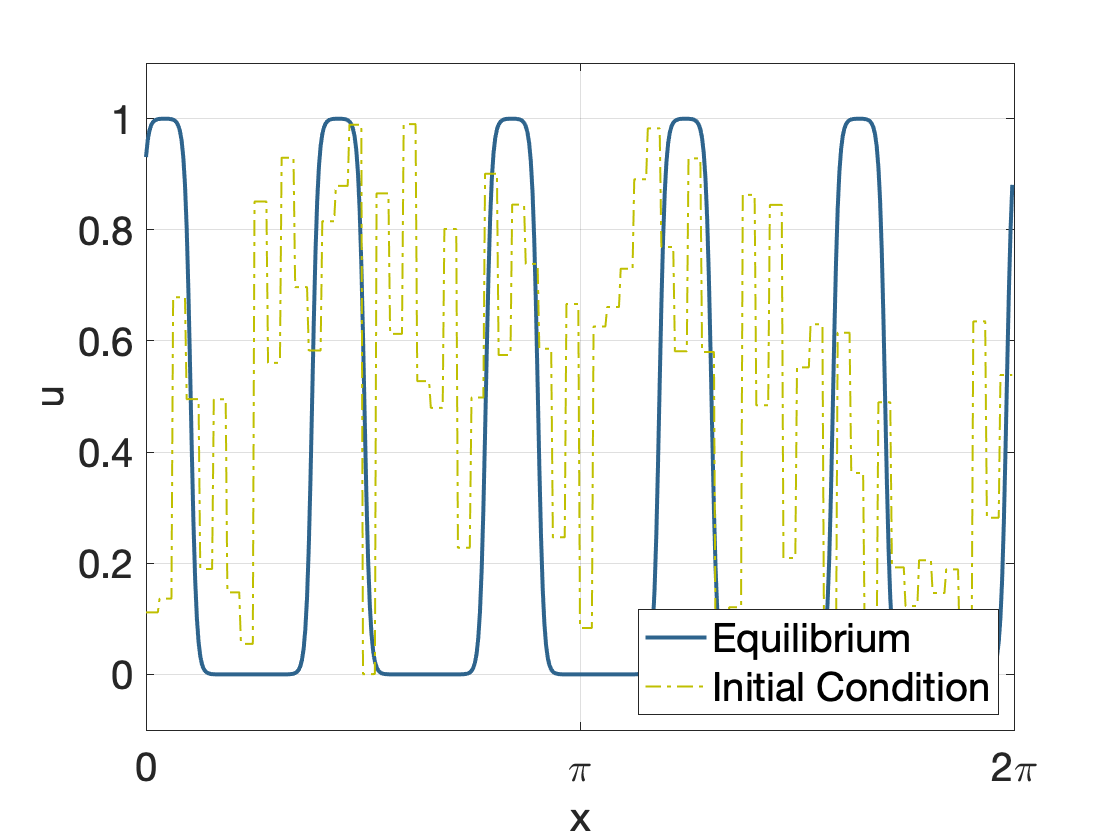}
\end{center}
\caption{Numerical evidence that $N$-bubble functions of equal size and equal distance are the minimizers of $E_0^{\text{NOK}}$ over $\mathcal{A}_N^{\omega}$. Left: $\mathcal{L}_{\delta}$ is the nonlocal operator with a constant kernel as shown in Example \ref{example:constantkernel}. Starting from an random initial (dashed dark green), the $L^2$ gradient flow dynamics leads us to an equilibrium (solid blue) with four bubbles of equal size and equal distance. In this example, $\omega = 0.2, \delta = 0.3$.  Right: $\mathcal{L}_{\delta}$ is the nonlocal operator with a power kernel ($\alpha = 2$) as shown in Example \ref{example:powerkernel}. The $L^2$ gradient flow dynamics results in equilibrium with five bubbles of equal size and equal distance. In this example, $\omega = 0.3, \delta = 0.3$. In both simulations, $\epsilon = 10\Delta x$ is fixed with $\Delta x = \frac{2\pi}{1024}$.}
\label{fig:numericalevidence} 
\end{figure}

\subsection{Minimizer of $E_0^{\emph{NOK}}$ over $\mathcal{A}_N^{\omega}$} \label{subsection:AOmega}

In this section, we consider the minimizer of $E_0^{\text{NOK}}$ over the periodic step functions in $\mathcal{A}_N^{\omega}$. Unlike the case $\mathcal{L}_{\delta} = -\Delta$, in which Green's function has analytical quadratic form and consequently one can prove that the minimizer over $\mathcal{A}_N^{\omega}$ is the $N$-bubble function of equal size and equal distance \cite{RenTruskinovsky_Elasticity2000, FrankLieb_LettersMathPhysics2019}, we are lack of analytical tools to characterize the minimizers of $E_0^{\text{NOK}}$ over $\mathcal{A}_N^{\omega}$.

On the other hand, starting from random initials, we numerically implement the $L^2$ gradient flow dynamics for the diffuse interface NOK model (\ref{functional:NOK}) with various nonlocal operators $\mathcal{L}_{\delta}$, and always obtain the $N$-bubble functions of equal size and equal distance as the equilibria. Figure \ref{fig:numericalevidence}  presents two numerical examples for the $L^2$ gradient flow dynamics. Here we use Fourier spectral method for spatial discretization and BDF2 for the temporal discretization. The two nonlocal operators are chosen as the one with constant kernel and the one with power kernel, see subsection \ref{subsubsection:preliminaries} for the preliminaries of the 1D nonlocal operators and examples \ref{example:powerkernel} and \ref{example:constantkernel} for power and constant kernels.  

The numerical results make us believe that $U_{N}^{\omega}$ is most probably the minimizer  of $E_0^{\text{NOK}}$ over $\mathcal{A}_N^{\omega}$. Thereby we propose the following conjecture:

\begin{conjecture}\label{theorem:equalsize_equaldistance}
Given $\omega \ll 1$ and $\gamma\gg 1$. $u^*(x)$ is the unique minimizer of $E_0^{\emph{NOK}}$ in $\mathcal{A}_N^{\omega}$ if and only if $u^*(x)$ is of equal size and equal distance. In other words, 
\begin{align}\label{eqn:conjucture}
U_N^{\omega} = \underset{u\in \mathcal{A}_N^{\omega}}{\mathrm{argmin}} E_0^{\emph{NOK}}[u] .
\end{align}
\end{conjecture}

The conjecture is provable for $N=2$ under the conditions that Green's function of $\mathcal{L}_{\delta}$ is symmetric $G_{\delta}(x,y) = g_{\delta}(|x-y|)$ and $g_{\delta}(\cdot)$ is strictly convex. However, these conditions can only be guaranteed when $\delta$ is close to 0. If $\delta$ is far from 0, say $\delta = 1$,  $g_{\delta}(\cdot)$ may lose the strict convexity. We will search for more analytical tools to prove the conjecture in our future work.

\subsection{Minimizer of $E_0^{\emph{NOK}}[U_N^{\omega}]$ over $N$} \label{subsection: UOmega}

Though it is mysterious to us whether the conjecture \ref{theorem:equalsize_equaldistance} holds theoretically, the numerical results convince us that it is still worth considering the second step of the minimization, minimizing $E_0^{\text{NOK}}[U_{N}^{\omega}]$ over $N\in\mathbf{Z}_+$. In this section, we will characterize the minimizers of of $E_0^{\text{NOK}}[U_N^{\omega}]$ over $N\in\mathbb{Z}^+$ for various nonlocal operators $\mathcal{L}_{\delta}$. Then we will study the influence of the nonlocal parameters $\gamma, \delta$ on the minimizers. To begin with, we briefly recall some preliminaries for the 1D nonlocal operators \cite{Du_Book2020, DuYang_SINA2016}.

\subsubsection{Preliminaries on the 1D nonlocal operators}\label{subsubsection:preliminaries}

In the 1D case, the nonlocal operator (\ref{kernel_Lij}) can be rewritten as
\begin{align}
\mathcal{L}_{\delta} u(x) = \int_{|s|\le \delta} \rho_{\delta}(s) (u(x+s)-u(x)) ds.
\end{align}
For the sake of simplicity, we consider the radially symmetric kernel
\begin{align}
\rho_{\delta}(s) = \frac{1}{\delta^3} \rho\left( \frac{|s|}{\delta}  \right), \quad \forall s\in[-\delta,\delta],
\end{align}
in which $\rho = \rho(\xi)$ is a nonnegative nonincreasing function with a compact support in $[0,1]$ and a bounded second moment,
\[
\int_0^1 \rho(\xi) \xi^2 d\xi = 1.
\]
Under periodic boundary conditions,  $e^{inx}$ is an eigenfunction of $\mathcal{L}_{\delta}$ with the corresponding eigenvalue
\begin{align}
\lambda_{\delta}(n) = \frac{2}{\delta^2} \int_0^1 \rho(\xi)(1-\cos n\delta\xi) d\xi = 2 \int_0^{\delta} \rho_{\delta}(s)(1-\cos ns) ds.
\end{align}

Here we list several examples of the nonlocal operators.

\begin{figure}[t] 
\begin{center}
\includegraphics[width=0.9\linewidth]{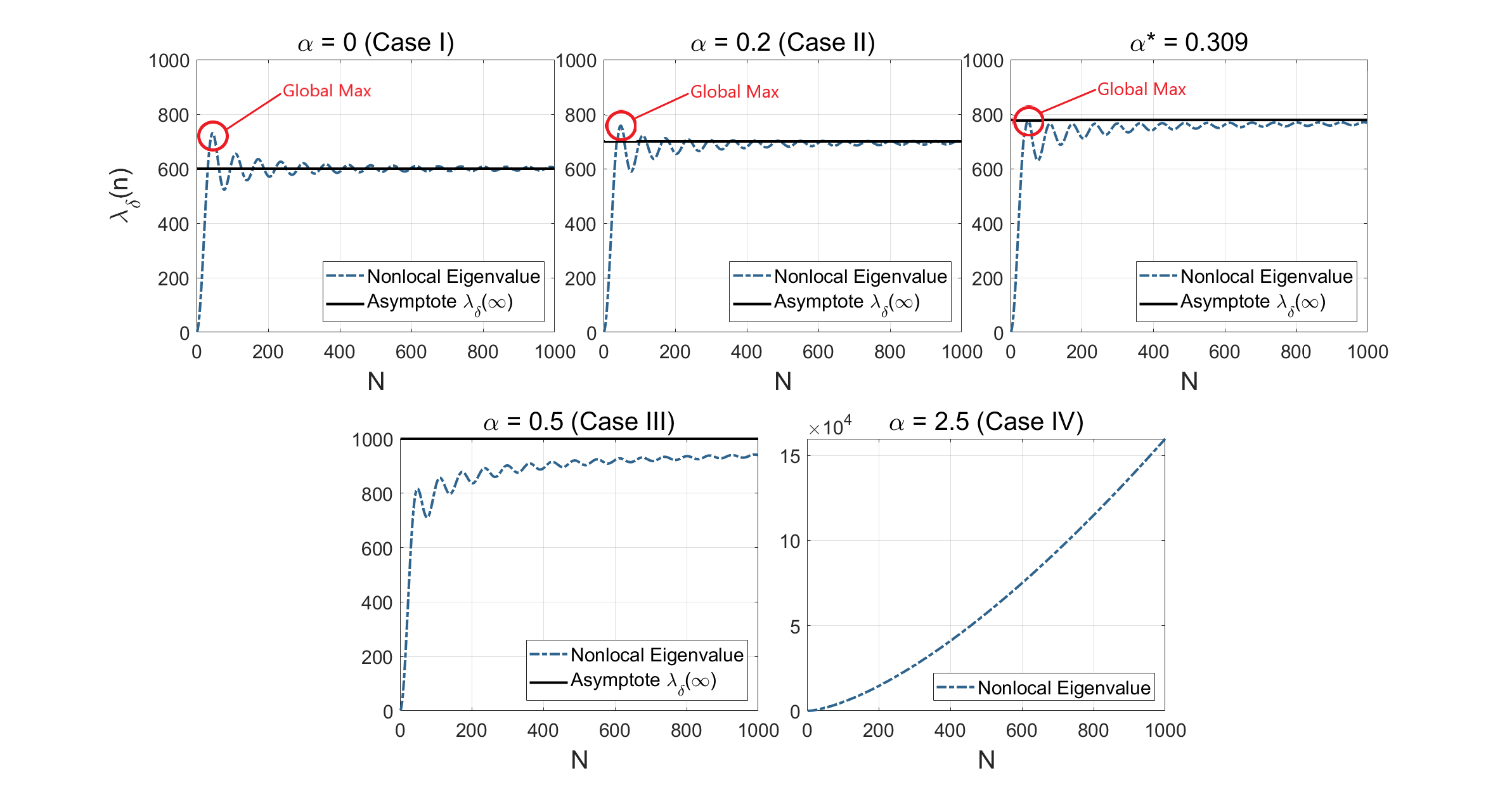}
\end{center}
\caption{Eigenvalues for power kernel operator $\rho_{\delta}(s) = \frac{3-\alpha}{\delta^{3-\alpha}|s|^{\alpha}}, s \in [-\delta,0)\cup(0,\delta]$.  $\alpha = 0$ is for Case I in which the eigenvalues oscillate about $\lambda_{\delta}(\infty)$. $\alpha = 0.2\in(0,\alpha^*)$ is a example for Case II in which eigenvalues initially oscillate around $\lambda_{\delta}(\infty)$ but eventually become oscillating from below and converging to $\lambda_{\delta}(\infty)$. $\alpha^*\approx 0.309$ is the dividing value of $\alpha$ between Case II and Case III. $\alpha = 0.5\in(\alpha^*,1)$ is a example for Case III in which the eigenvalues oscillate from below $\lambda_{\delta}(\infty)$ and converge to $\lambda_{\delta}(\infty)$. $\alpha = 2.5\in[1,3)$ is a example for Case IV in which eigenvalues monotonically grow to $\infty$. The red circles indicate whether $\lambda_{\delta}(n)$ attains a global maximum, which is an important feature to determine whether the optimal solution $N^*(\gamma)$ of (\ref{eqn:min_N}) has an upper bound.  }
\label{fig:eigenvalue_powerkernel} 
\end{figure}

\begin{example}[Nonlocal operator with power kernel] \label{example:powerkernel}
The nonlocal operator $\mathcal{L}_{\delta}$ with the power kernel
\begin{align}\label{eqn:Kernel_Power}
\rho(\xi) = \frac{3-\alpha}{\xi^{\alpha}},\ \xi\in(0,1), {\text{\ or}} \quad \rho_{\delta}(s) = \frac{3-\alpha}{\delta^{3-\alpha}|s|^{\alpha}}, \ s \in [-\delta,0)\cup(0,\delta], \quad \text{where\ } \alpha \in (0,3),
\end{align}
has eigenvalues
\begin{align}\label{eqn:eigenvalue}
\lambda_{\delta}(n)=2(3-\alpha)\delta^{-(3-\alpha)}|n|^{\alpha-1}\int_{0}^{n\delta}\frac{1-\cos{t}}{t^{\alpha}}dt.
\end{align}
When $\alpha\in(0,1)$, the eigenvalue $\lambda_{\delta}(n)$ has a finite limit
\[
\lambda_{\delta}(\infty):=\lim_{n\to\infty}\lambda_{\delta}(n)=\frac{2(3-\alpha)}{\delta^2(1-\alpha)}.
\]

In figure \ref{fig:eigenvalue_powerkernel}, we plot $\lambda_{\delta}(n)$ for several values of $\alpha =0, 0.2, 0.309, 0.5, 2.5$ with a fixed $\delta=0.1$.
\end{example}

\begin{example}[Nonlocal operator with constant kernel] \label{example:constantkernel}
The nonlocal operator $\mathcal{L}_{\delta}$ with the constant kernel
\begin{align}\label{eqn:Kernel_Const}
\rho(\xi) \equiv 3, \ {\text{\ or}} \quad \rho_{\delta}(s) = \frac{3}{\delta^3},
\end{align}
has eigenvalues
\begin{align}\label{eqn:Eigenvalue_Const}
\lambda_{\delta}(n) = \frac{6}{\delta^2}\left(1-\emph{sinc} (n\delta) \right).
\end{align}
with a finite limit
\[
\lambda_{\delta}(\infty):=\lim_{n\to\infty}\lambda_{\delta}(n) = \frac{6}{\delta^2}.
\]
The constant kernel can be regarded as the degenerate power kernel when $\alpha \rightarrow 0$.
\end{example}

\begin{example}[Classification of the nonlocal operator with power kernel for $\alpha\in[0,3)$]\label{example:classification}
 For nonlocal kernels (\ref{eqn:Kernel_Power}) and (\ref{eqn:Kernel_Const}), we can classify them into four cases according to the behavior of the eigenvalues $\lambda_{\delta}(n)$.
\begin{itemize}
\item[Case \emph{I}.]  $\alpha = 0$. In this case, the eigenvalues $\lambda_{\delta}(n)$ oscillate around the constant value $\lambda_{\delta}(\infty)$ and converge to $\lambda_{\delta}(\infty)$ as $n\to\infty$. 

\item[Case \emph{II}.]  $\alpha \in (0,\alpha^*)$. In this case, the eigenvalues $\lambda_{\delta}(n)$ initially oscillate around the constant value $\lambda_{\delta}(\infty)$, but eventually oscillate from below and converge to $\lambda_{\delta}(\infty)$. Here the value of $\alpha^*\in(0.308,0.309)$ will be justified by the lemma \ref{lemma:alpha*} below. 

\item[Case \emph{III}.] $\alpha \in (\alpha^*, 1)$. In this case, the eigenvalues $\lambda_{\delta}(n)$ oscillate from below and converge to the limit $\lambda_{\delta}(\infty)$  as $n\to\infty$. 
\item[Case \emph{IV}.] $\alpha \in (1, 3)$. In this case, the eigenvalues $\lambda_{\delta}(n)$ monotonically increase to $\infty$ as $n \rightarrow \infty$. 
\end{itemize}

\end{example}

\begin{example}[Nonlocal operator with Gauss kernel] The nonlocal operator $\mathcal{L}_{\delta}$ with the Gauss-type kernel
\begin{align}\label{eqn:Kernel_Gauss}
\rho(\xi) = \frac{2}{\sqrt{\pi}}e^{-\xi^2}, \ {\text{\ or}} \quad \rho_{\delta}(s) = \frac{2}{\sqrt{\pi}\delta^3}e^{-\frac{|s|^2}{\delta^2}},
\end{align}
has eigenvalues
\begin{align}\label{eqn:Eigenvalue_Gauss}
\lambda_{\delta}(n) = \frac{4}{\delta^2}\left( 1 - e^{-\frac{(n\delta)^2}{4}} \right).
\end{align}
Because of  its exponential decay as $|s|\rightarrow\infty$, $\rho_{\delta}$ can be regarded as being with a compact support in the finite domain $\mathbb{T}^d$ when the domain size is much larger than $\delta$, and consequently can be periodically extended.

\end{example}

We conclude the preliminaries by introducing the following lemma to determine the critical value $\alpha^*$ for the power kernel.

\begin{lemma}\label{lemma:alpha*}
Let the nonlocal operator $\mathcal{L}_{\delta}$ be with the power kernel (\ref{eqn:Kernel_Power}) and $\alpha\in(0,1)$. Treating $n$ as a continuum, the eigenvalues (\ref{eqn:eigenvalue}) are oscillatory, and  there exists a critical value $\alpha^* \in (0.308, 0.309)$ such that for any $\alpha \in (0,\alpha^*)$, the global maximum of $\lambda_{\delta}(n)$ is reached at a finite value $N^*<\infty$, and for any $\alpha\in(\alpha^*,1)$, the supremium of $\lambda_{\delta}(n)$ is reached at $n\rightarrow\infty$.
\end{lemma}

\begin{proof}
The eigenvalues of the nonlocal operator $\mathcal{L}_{\delta}$ with power kernel can be calculated as 
\begin{align*}
    \lambda_{\delta}(n) & = 2(3-\alpha)\delta^{\alpha-3}n^{\alpha-1}\int_{0}^{n\delta}\frac{1-\cos{t}}{t^{\alpha}}dt \\
    & = 2(3-\alpha)\delta^{\alpha-3}n^{\alpha-1}\left(\int_{0}^{n\delta}\frac{1}{t^{\alpha}}dt - \int_{0}^{n\delta}\frac{\cos{t}}{t^{\alpha}}dt\right) \\
    & = \frac{2(3-\alpha)}{\delta^2(1-\alpha)} - 2(3-\alpha)\delta^{\alpha-3}n^{\alpha-1}\int_{0}^{n\delta}\frac{\cos{t}}{t^{\alpha}}dt.
\end{align*}
The oscillation of $\lambda_{\delta}(n)$ is obviously determined by the oscillation of the integral term $\int_{0}^{n\delta}\frac{\cos{t}}{t^{\alpha}}dt$ as $n\rightarrow\infty$.
Note that for $\alpha \in (0,1)$, $\int_{0}^{\infty}\frac{\cos{t}}{t^{\alpha}}dt$ is convergent, therefore 
$
   \lambda_{\delta}(\infty) = \frac{2(3-\alpha)}{\delta^2(1-\alpha)} . 
$

To find the critical value $\alpha^*$ such that for any $\alpha\in(\alpha^*,1)$, $\lambda_{\delta}(n)$ reaches the supremium at $\infty$, one needs for any $n>0$ that $\lambda_{\delta}(n) < \lambda_{\delta}(\infty)$, or equivalently $\int_{0}^{n\delta}\frac{\cos{t}}{t^{\alpha}}dt > 0$. Defining $f(x) = \int_{0}^{x}\frac{\cos{t}}{t^{\alpha}}dt$, simple calculation yields that $f(x)$ attains the global minimum at $x = \frac{3\pi}{2}$. Consequently, the critical value $\alpha^*$ is determined by the condition $f(\frac{3\pi}{2}) = 0$. Note that $f(\frac{3\pi}{2})= \int_{0}^{\frac{3\pi}{2}}\frac{\cos{t}}{t^{\alpha}}dt$ is an increasing function with respect to $\alpha\in(0,1)$, and
\[
\lim_{\alpha\rightarrow0}\int_{0}^{\frac{3\pi}{2}}\frac{\cos{t}}{t^{\alpha}} dt = -1, \quad \lim_{\alpha\rightarrow1}\int_{0}^{\frac{3\pi}{2}}\frac{\cos{t}}{t^{\alpha}} dt = \infty,
\]
there exists a unique $\alpha^*$ such that $f(\frac{3\pi}{2}) = 0$. Using bisection method, we numerically find that $\alpha^* \in (0.308, 0.309)$.    

Besides, when $\alpha \in (0,\alpha^*)$, the integral $\int_{0}^{\infty}\frac{\cos{t}}{t^{\alpha}}dt$ changes the signs for small value of $n$, but eventually becomes positive. Therefore the eigenvalues $\lambda_{\delta}(n)$ are initially oscillatory about $\lambda_{\delta}(\infty)$ but eventually become smaller than and asymptotically approach $\lambda_{\delta}(\infty)$. 
\end{proof}

\subsubsection{Formulation of $E_0^{\emph{NOK}}[U_N^{\omega}]$}

For a fixed $N$, the NOK energy $E_0^{\text{NOK}}[U_N^{\omega}]$ for the equal-area and equal-distance bubble function has a simple form as given in the following lemma.

\begin{lemma}\label{lemma:energy_U}
For a fixed $N$, the NOK energy $E_0^{\emph{NOK}}[U_N^{\omega}]$ can be represented as a Fourier series, involving attractive term $E^{\emph{att}}$ and repulsive term $E^{\emph{rep}}$,
\begin{align}\label{eqn:energy_Fourierseries}
E_0^{\emph{NOK}}[U_N^{\omega}] = E^{\emph{tot}}: = E^{\emph{att}} + E^{\emph{rep}}= 2N + \frac{2\gamma}{\pi}\sum_{m=1}^{\infty}\frac{1}{\lambda_{\delta}(mN)}\left( \frac{\sin(m\pi\omega)}{m} \right)^2.
\end{align}
\end{lemma}
\begin{proof}
The first term $2N$ is fixed. It remains to work on the second term $E^{\text{rep}}$ in (\ref{eqn:energy_Fourierseries}). Consider the Fourier series of any function $u\in\mathcal{A}_N^{\omega}$:
\[
u(x) = \sum_{n=-\infty}^{\infty} \hat{u}_n e^{-inx}.
\]
Since the eigenvalue of the nonlocal operator $\mathcal{L}_{\delta}$ is $\lambda_{\delta}(n)$, the Fourier series of $\mathcal{L}_{\delta}^{-1}u$ is given as
\[
(\mathcal{L}_{\delta}^{-1}u)(x) = \sum_{n=-\infty}^{\infty} \frac{1}{\lambda_{\delta}(n)}\hat{u}_n e^{-inx},
\]
and then by (\ref{eqn:SharpLimit_1d})
\begin{align}\label{eqn:E_2}
E^{\text{rep}} =  \dfrac{\gamma}{2}\int_0^{2\pi} (\mathcal{L}_{\delta})^{-1}(u(x)-\omega) (u(x)-\omega)\ \text{d}x 
=  2\gamma \pi \sum_{n=1}^{\infty} \frac{1}{\lambda_{\delta}(n)} |\hat{u}_n|^2.
\end{align}
Now taking $u(x) = U_N^{\omega}(x)$, then the Fourier coefficient $\hat{u}_n$ becomes
\begin{align*}
\hat{u}_n  & = \frac{1}{2\pi}\int_{0}^{2\pi}u(x)e^{inx}dx \\
              & =\frac{1}{2\pi} \sum_{j=0}^{N-1}\left(\int_{\frac{2\pi j}{N}}^{\frac{2\pi j}{N}+\frac{2\pi\omega}{N}}u(x)e^{inx}dx+\int_{\frac{2\pi j}{N}+\frac{2\pi\omega}{N}}^{\frac{2\pi (j+1)}{N}}u(x)e^{inx}dx\right) \\
              & = \frac{1}{2\pi} \sum_{j=0}^{N-1}\int_{\frac{2\pi j}{N}}^{\frac{2\pi j}{N}+\frac{2\pi\omega}{N}}e^{inx}dx \\
              & = \frac{1}{2\pi} \sum_{j=0}^{N-1}\frac{1}{in} e^{in\frac{2\pi j}{N}}\left[e^{in\frac{2\pi\omega}{N}}-1\right] \\
              & = \frac{1}{2\pi} \frac{1}{in}\frac{1-e^{in2\pi}}{1-e^{in\frac{2\pi}{N}}}\left[e^{in\frac{2\pi\omega}{N}}-1\right].
\end{align*}
Since $e^{in(2\pi)}=1$, so we have $1-e^{in(2\pi)}=0$, $\forall n$. If $\frac{n}{N}\not\in\mathbb{Z}$, then $1-e^{in\frac{2\pi}{N}}\neq0$, so that $\hat{u}_n = 0$. If $\frac{n}{N}\in\mathbb{Z}$, then $1-e^{in\frac{2\pi}{N}} = 0$, and $\frac{1-e^{in2\pi}}{1-e^{in\frac{2\pi}{N}}} = N$, resulting in $\hat{u}_n = \frac{1}{2\pi}\frac{1}{in}N\left[e^{in\frac{2\pi\omega}{N}}-1\right]$. Therefore we have
\begin{align}\label{eqn:fourier_coefficient_square}
 |\hat{u}_n|^2=\begin{cases}
     \frac{1}{\pi^2}\frac{N^2}{n^2}\cdot \sin^2{\frac{n\pi\omega}{N}},& \frac{n}{N}\in\mathbb{Z}.\\
     0, &  \frac{n}{N}\not\in\mathbb{Z}.
 \end{cases}
\end{align}
Denoting $\frac{n}{N} = m$ and inserting $|\hat{u}_n|^2$ into $E^{\text{rep}} $ in (\ref{eqn:E_2}) yields the Fourier series (\ref{eqn:energy_Fourierseries}).
\end{proof}

According to lemma \ref{lemma:energy_U}, the minimizer of
\begin{align}\label{eqn:min_N}
\min_{N\in\mathbb{Z}^+} E^{\text{tot}} = E^{\text{att}} + E^{\text{rep}} = 2N + \frac{2\gamma}{\pi}\sum_{m=1}^{\infty}\frac{1}{\lambda_{\delta}(mN)}\left(\frac{\sin(m\pi\omega)}{m}\right)^2,
\end{align}
is determined by the competition between $E^{\text{att}} = 2N$, which is a linear term of $N$, and $E^{\text{rep}}$, whose growth with respect to $N$ is controlled by the behavior of the eigenvalues $\lambda_{\delta}(n)$ of the nonlocal operator $\mathcal{L}_{\delta}$.

In what follows, we will study the nonlocal effect on the minimizers. More specifically, we will consider the four cases listed in the example \ref{example:classification}. For the sake of simplicity, we denote 
\begin{align}\label{eqn:F}
F(N;\delta,\omega) := \sum_{m=1}^{\infty}\frac{1}{\lambda_{\delta}(mN)}\left(\frac{\sin(m\pi\omega)}{m}\right)^2,
\end{align}
such that $E^{\text{rep}} = \frac{2\gamma}{\pi}F(N;\delta,\omega)$.

\subsection{Nonlocal $\gamma$-Effect on Optimizers for Power Kernel}\label{subsection:PowerKernel}

In this section, we focus on the nonlocal operator with power kernel as discussed in Example \ref{example:classification}, and explore its nonlocal effect on the pattern of minimizers for equation (\ref{eqn:min_N}). We reveal that the optimal number $N^*(\gamma)$ of bubbles for the minimization problem  (\ref{eqn:min_N}) is always a nondecreasing function with respect to the long-range repulsion strength $\gamma$, and more importantly, $N^*(\gamma)$ may have an upper bound under the condition that $\lambda_{\delta}(n)$ attains a global minimum at some finite $n>0$ (Case I in Example \ref{example:classification}). On the other hand, if $\lambda_{\delta}(n)$ does not attain a global minimum at any finite $n>0$ (Cases III \& IV in Example \ref{example:classification}), $N^*(\gamma)$ becomes unbounded from above. We discuss the four cases in Example \ref{example:classification} in the following subsections.

\subsubsection{$\gamma$-Effect on Minimizers for Case I}

We begin with the nonlocal operator $\mathcal{L}_{\delta}$ with a constant kernel as illustrated in Example \ref{example:constantkernel}, which is Case I in Example \ref{example:classification}. We will perform analysis on the optimal number $N^*$ of bubbles for the minimization problem (\ref{eqn:min_N}), and the dependence of $N^*$ on the repulsion strength $\gamma$. We reveal that there exists an upper bound for $N^*$ as $\gamma \rightarrow \infty$:
\begin{align}
N^*(\gamma) \le \bar{N}, \quad \text{for any\ } \gamma \in \mathbb{R}^+ .
\end{align}

Given the kernel function (\ref{eqn:Kernel_Const}), the eigenvalues of $\mathcal{L}_{\delta}$ associated with this kernel is given by \cite{DuYang_SINA2016} as shown in (\ref{eqn:Eigenvalue_Const}). Inserting the eigenvalue $\lambda_{\delta}$  into $F(N;\delta,\omega) $, we have
\begin{align*}
F(N;\delta,\omega) 
=\frac{\delta^2}{6}\sum_{m=1}^{\infty} \frac{1}{1-\text{sinc}{(mN\delta)}}\left(\frac{\sin{m\pi\omega}}{m}\right)^2.
\end{align*}
We denote by $F^{\infty} := F(\infty; \delta,\omega)$ the limit of $F(N;\delta,\omega)$ as $N\rightarrow\infty$. Then it is evident that
\begin{align}
F^{\infty}=\frac{\delta^2}{6}\sum_{m=1}^{\infty}\left(\frac{\sin{m\pi\omega}}{m}\right)^2>0.
\end{align}

The optimal number $N^*(\gamma)$ for the minimization problem (\ref{eqn:min_N}) is closely related to $F(N;\delta,\omega)$ in \ref{eqn:F}. We will show that $F(N;\delta,\omega)$ behaves like the sinc function. To this end, we need several lemmas regarding the properties of $F(N;\delta,\omega)$.

\begin{lemma}\label{lemma:E_infty}
Let $N\in[1,\infty)$ be a continuum, and $\omega\in(0,1/2]$.  For $\forall k\in\mathbb{Z}_+$, $F\left(N;\delta,\omega\right)|_{N=\frac{k\pi}{\delta}}=F^{\infty}$. 
\end{lemma}
\begin{proof}
When $N\delta=k\pi$, the Fourier series $F(N;\delta,\omega)$ becomes
\begin{align*}
F(N;\delta,\omega)  
&=\frac{\delta^2}{6}\sum_{m=1}^{\infty} \frac{mN\delta}{mN\delta-\sin{mN\delta}}\left(\frac{\sin{m\pi\omega}}{m}\right)^2 \\
&=\frac{\delta^2}{6}\sum_{m=1}^{\infty}\frac{mk\pi}{mk\pi-\sin{mk\pi}}\left(\frac{\sin{m\pi\omega}}{m}\right)^2
=\frac{\delta^2}{6}\sum_{m=1}^{\infty}\left(\frac{\sin{m\pi\omega}}{m}\right)^2
=F^{\infty},
\end{align*}
completing the proof.
\end{proof}

\begin{lemma}\label{lemma:E_Oscilate}
Let $N\in[1,\infty)$ be a continuum, and $\omega\in(0,1/2]$. $F(N;\delta,\omega)$ decreases at $N=\frac{(2k-1)\pi}{\delta}$ and increases at $N = \frac{2k\pi}{\delta}$ for any $\forall k\in \mathbb{Z}_+$.
\end{lemma}
\begin{proof}
Since $N\in[1,\infty)$ is treated as a continuum, one can take the derivative of $F(M;\delta,\omega)$ with respect to $N$,
\begin{align*}
\frac{\partial F(N;\delta,\omega)}{\partial N} = \frac{\delta^2}{6}\sum_{m=1}^{\infty}\frac{mN\delta\cos{mN\delta}-\sin{mN\delta}}{(mN\delta-\sin{mN\delta})^2}(m\delta)\left(\frac{\sin{m\pi\omega}}{m}\right)^2.
\end{align*}
If $N\delta=2k\pi$, then $\cos{mN\delta}=1$, $\sin{mN\delta}=0$ for any $m\in\mathbb{Z}^+$, so we have: 
\begin{align*}
\frac{\partial F(N;\delta,\omega)}{\partial N}\Big|_{N\delta=2k\pi}  =\frac{\delta^2}{6}\sum_{m=1}^{\infty}\frac{mN\delta}{(mN\delta)^2}(m\delta)\left(\frac{\sin{m\pi\omega}}{m}\right)^2 =\frac{\delta^2}{6N}\sum_{m=1}^{\infty}\left(\frac{\sin{m\pi\omega}}{m}\right)^2>0,
\end{align*}
Hence, $F(M;\delta,\omega)$ increases at $N=\frac{2k\pi}{\delta}$. If $N\delta=(2k+1)\pi$, then $\cos{mN\delta}=(-1)^m$, $\sin{mN\delta}=0$. Noting the following trigonometric identity,
\[
\sum_{n=1}^{\infty}(-1)^n\left(\frac{\sin{nx}}{n}\right)^2=-\frac{x^2}{2}, \quad x\in(-\pi,\pi),
\]
we have
\begin{align*}
\frac{\partial F(N;\delta,\omega)}{\partial N}\Big|_{N\delta=2k\pi+1}
=\frac{\delta^2}{6N}\sum_{m=1}^{\infty}(-1)^m\left(\frac{\sin{m\pi\omega}}{m}\right)^2
=-\frac{\delta^2}{6N}\frac{(\pi\omega)^2}{2}<0,
\end{align*}
which implies that $F(N;\delta,\omega)$ decreases at $N = \frac{(2k+1)\pi}{\delta}$.
\end{proof}

\begin{lemma}\label{lemma:E_OppositeSign}
Let $N\in[1,\infty)$ be a continuum, and $\omega\in(0,1/2]$ and $\delta\in (0,\pi]$.  We have that $F(N;\delta,\omega)>E^{\infty}$ for $N\in\left(\frac{2k\pi}{\delta},\frac{(2k+1)\pi}{\delta}\right)$, and $F(N;\delta,\omega)<E^{\infty}$ for $N\in \left(\frac{(2k+1)\pi}{\delta},\frac{(2k+2)\pi}{\delta} \right)$ with $ k=0,1,\cdots$.
\end{lemma}
\begin{proof}
Let $x = N\delta$ and $y = \pi\omega$. We have
\begin{align*}
F(N;\delta,\omega)-E^{\infty} = \frac{\delta^2}{6}\sum_{m=1}^{\infty}\frac{\sin{mx}}{mx-\sin{mx}}\left(\frac{\sin{my}}{m}\right)^2.
\end{align*}
According to lemma \ref{lemma:E_infty}, we know that $E^{\text{rep}}(N;\delta,\omega)-E^{\infty} = 0$ at $x = k\pi$. Furthermore, for $x = k\pi \pm \epsilon$ with $\epsilon\in(0,\pi)$, we have
\begin{align*}
   \Big( F(N;\delta,\omega)-E^{\infty} \Big) \Big|_{x = k\pi-\epsilon} &  = - \frac{\delta^2}{6}\sum_{m=1}^{\infty}\frac{(-1)^{mk}\sin{m\epsilon}}{mk\pi  -  (m\epsilon - (-1)^{mk}\sin{m\epsilon})}\left(\frac{\sin{my}}{m}\right)^2, \\
   \Big( F(N;\delta,\omega)-E^{\infty} \Big) \Big|_{x = k\pi+\epsilon} &  = + \frac{\delta^2}{6}\sum_{m=1}^{\infty}\frac{(-1)^{mk}\sin{m\epsilon}}{mk\pi  +  (m\epsilon - (-1)^{mk}\sin{m\epsilon})}\left(\frac{\sin{my}}{m}\right)^2.
\end{align*}
Consequently, $F(N;\delta,\omega)-E^{\infty}$ has opposite sign for $x = k\pi \pm \epsilon$ with $\epsilon\in(0,\pi)$. What is more, since 
\[
mk\pi  -  (m\epsilon - (-1)^{mk}\sin{m\epsilon} ) < mk\pi  +  (m\epsilon - (-1)^{mk}\sin{m\epsilon}),
\]
the magnitude of $F(N;\delta,\omega)-E^{\infty}$ at $x = k\pi-\epsilon$ is greater than that at $x = k\pi+\epsilon$. 

It remains to show that $F(N;\delta,\omega)-E^{\infty} > 0$ for $x\in(0,\pi)$, then the results can be guaranteed by the periodic sign change of $F(N;\delta,\omega)-E^{\infty}$. Consider the power series representation of $F(N;\delta,\omega)-E^{\infty}$, we have that 
\begin{align*}
F(N;\delta,\omega)-E^{\infty} = \frac{\delta^2}{6}\sum_{m=1}^{\infty}\left[\frac{\sin{mx}}{mx}+\left(\frac{\sin{mx}}{mx}\right)^2+\left(\frac{\sin{mx}}{mx}\right)^3+\cdots\right]\left(\frac{\sin{my}}{m}\right)^2
\end{align*}
Define 
\begin{align*}
    g_n(x) = \frac{\delta^2}{6}\sum_{m=1}^{\infty}\left(\frac{\sin{mx}}{mx}\right)^n\left(\frac{\sin{my}}{m}\right)^2,
\end{align*}
we have 
\begin{align*}
g_{2k}(x)+g_{2k+1}(x)  = \frac{\delta^2}{6}\sum_{m=1}^{\infty}\left(\frac{\sin{mx}}{mx}\right)^{2k}\left(\frac{\sin{my}}{m}\right)^2 \left( 1 + \frac{\sin{mx}}{mx} \right) \ge 0.
\end{align*}
For $g_1(x) = \frac{\delta^2}{6}\sum_{m=1}^{\infty}\left(\frac{\sin{my}}{m}\right)^2\frac{\sin{mx}}{mx}$, we can verify that for $x\in(0,\pi)$ and $y\in (0,\frac{\pi}{2}]$, the piecewise function
\begin{align*}
f(x) = 
\begin{cases}
\frac{y(\pi-y)}{2}x - \frac{\pi}{8}x^2, \quad x\in(0, 2y] , \\
\frac{y^2}{2}(\pi-x), \hspace{0.44in} x\in(2y,\pi),
\end{cases}
\end{align*}
is associated with a convergent Fourier series
\[
f(x) = \sum_{m=1}\frac{\sin^2{my}}{m^3}\sin{mx}.
\]
Hence, for $x\in(0,\pi)$, $y\in(0,\frac{\pi}{2}]$,
\begin{align*}
    g_1(x) =  \frac{\delta^2}{6}\sum_{m=1}^{\infty}\left(\frac{\sin{my}}{m}\right)^2\frac{\sin{mx}}{mx} = \frac{\delta^2}{6x}f(x) > 0.
\end{align*}
Therefore, we conclude that for $x\in(0,\pi)$, 
\begin{align*}
F(N;\delta,\omega)-E^{\infty} = g_1(x) + (g_2(x)+g_3(x)) + (g_4(x)+g_5(x)) + \cdots \geq g_1(x) > 0.
\end{align*}
The desired result is obtained.
\end{proof}

\begin{lemma}
Let $N\in[1,\infty)$ be a continuum, $\omega\in(0,1/2]$ and $\delta\in(0,\pi]$. $F(N;\delta,\omega)$ achieves the global minimum at some finite $\tilde{N}<\infty$.
\end{lemma}
\begin{proof}
According to lemmas \ref{lemma:E_infty}-\ref{lemma:E_OppositeSign}, we have
\begin{align*}
F(N;\delta,\omega) 
\begin{cases}
> E^{\infty}, \quad & N\delta \in (2k\pi, 2k\pi+\pi), \\
< E^{\infty}, \quad & N\delta \in (2k\pi+\pi, 2k\pi+2\pi), \\
= E^{\infty}, \quad & N\delta = k\pi,
\end{cases}
\end{align*}
and the magnitude of $F(N;\delta,\omega) - E^{\infty}$ over interval $(k\pi, (k+1)\pi)$ decreases to zero as $k\rightarrow\infty$. In other words, $F(N;\delta,\omega) - E^{\infty}$ behaves similarly as the sinc function.
Therefore $F(N;\delta,\omega)$ attains the global minimum at some point $\tilde{N}$ in the first trough.
\end{proof}

\begin{remark}
Since $F(N;\delta,\omega)$ behaves similarly as $\emph{sinc}(x)$, $F(N;\delta,\omega)$ reaches the global minimum at some finite value not only for the continuum variable $N$, but also for the discrete variable $N\in\mathbb{Z}^+$.
\end{remark}

Hereafter in this paper, we denote by $\tilde{N}$ the smallest optimal value of $N\in\mathbb{Z}^+$ for which $F(N;\delta,\omega)$ reaches the global minimum, and by $N^*(\gamma)$ the largest optimal value of $N\in\mathbb{Z}^+$ for which $E^{\text{tot}}$ reaches the global minimum.  The following theorem provides an upper bound for $N^*(\gamma)$.

\begin{theorem}\label{theorem:Case1_bound}
Let $N\in\mathbb{Z}^+$, $\omega\in(0,1/2)$ and $\delta\in(0,\pi]$. For the minimization problem (\ref{eqn:min_N}) in which $\lambda_{\delta}(n)$ is the eigenvalue of $\mathcal{L}_{\delta}$ with constant kernel (\ref{eqn:Kernel_Const}), the optimal number of bubbles $N^*(\gamma)$, as a function of $\gamma$, is a nondecreasing function with an upper bound
\begin{align}
N^*(\gamma) \le \tilde{N}, \quad \text{for\ } \forall \gamma >0.
\end{align}
\end{theorem}
\begin{proof}
Firstly, we prove the existence of the upper bound.  For the sake of brevity, we explicitly indicate the dependence of $E^{\text{tot}}$ on $N$ and $\gamma$, namely, $E^{\text{tot}} = E^{\text{tot}}(N,\gamma)$. Note that $E^{\text{att}} = 2N$ is linearly increasing with respect to $N$, and $E^{\text{rep}} = \frac{2\gamma}{\pi}F(N;\delta,\omega)$ reaches the global minimum at $\tilde{N}$, then it is evident that $E^{\text{tot}}(\tilde{N},\gamma) < E^{\text{tot}}(N,\gamma)$ for any $N>\tilde{N}$. Therefore $N^*(\gamma)$ does not exceed $\tilde{N}$.

Secondly, we show that $N^*(\gamma)$ is a nondecreasing function. we will prove that for any $\gamma_2 > \gamma_1$, 
\begin{align}
E^{\text{tot}}(N_2, \gamma_2) > E^{\text{tot}}(N^*(\gamma_1), \gamma_2), \quad \text{for\ } N_2 \in [1,N^*(\gamma_1)) \cap \mathbb{Z}^+,
\end{align}
which implies that $\min_{N\in\mathbb{Z}^+}E^{\text{tot}}(N, \gamma_2)$ cannot occur over $[1,N^*(\gamma_1))\cap \mathbb{Z}^+$.  Assuming for some $\gamma_2 > \gamma_1$, there exists a $N_2\in [1,N^*(\gamma_1))\cap \mathbb{Z}^+$ such that $E^{\text{tot}}(N_2, \gamma_2) \le E^{\text{tot}}(N^*(\gamma_1), \gamma_2)$, namely,
\[
2N_2+\frac{2\gamma_2}{\pi} F(N_2;\delta,\omega) \le 
2N^*(\gamma_1)+\frac{2\gamma_2}{\pi} F(N^*(\gamma_1);\delta,\omega).
\]
Then we have
\[
\frac{1}{\gamma_1} > \frac{1}{\gamma_2} \ge \frac{F(N_2;\delta,\omega) - F(N^*(\gamma_1);\delta,\omega)}{\pi(N^*(\gamma_1)-N_2)},
\]
and then
\[
2N_2+\frac{2\gamma_1}{\pi} F(N_2;\delta,\omega) < 
2N^*(\gamma_1)+\frac{2\gamma_1}{\pi} F(N^*(\gamma_1);\delta,\omega),
\]
which leads to $E^{\text{tot}}(N_2, \gamma_1) < E^{\text{tot}}(N^*(\gamma_1), \gamma_1)$, contradicting with the optimality of $N^*(\gamma_1)$. Therefore, the minimum $N^*(\gamma_2)$ for $\min_{N\in \mathbb{Z}^+}E^{\text{tot}}(N, \gamma_2)$ must be achieved over $[N^*(\gamma_1), \tilde{N}]\cap \mathbb{Z}^+$, resulting in the monotone increment of $N^*(\gamma)$.
\end{proof}

\subsubsection{$\gamma$-Effect on Minimizers for Case III} 

In the current and next subsections, we will skip over the discussion of case II and move towards cases III and IV for now as the results of cases III and IV will inspire and ease the discussion of case II. 

In case III, the eigenvalues $\lambda_{\delta}(n)$ oscillates from below and asymptotically approach a limit as $n\rightarrow \infty$. In this case, we have that
\begin{align}\label{eqn:infty_energy}
    F^{\infty} := \lim_{N\to\infty}F(N;\delta,\omega)=\sum_{m=1}^{\infty} \frac{1}{\lambda_{\delta}(\infty)} \left(\frac{\sin{m\pi\omega}}{m} \right)^2= \frac{\delta^2(1-\alpha)}{2(3-\alpha)}\sum_{m=1}^{\infty} \left(\frac{\sin{m\pi\omega}}{m} \right)^2.
\end{align}

\begin{lemma}
$F(N;\delta,\omega)$ does not reach a global minimum at any finite $N\in\mathbb{Z}^+$.
\end{lemma}
\begin{proof}
Note that $\lambda_{\delta}(n)$ is oscillating and asymptotically approach the limit $\lambda_{\delta}(\infty)$ from below, we have 
\[
\lambda_{\delta}(\infty) > \lambda_{\delta}(n), \quad n\in\mathbb{Z}^+.
\]
Then for any $N\in\mathbb{Z}^+$, 
\begin{align}\label{eqn:sub_energy}
F^{\infty} - F(N;\delta,\omega)=\sum_{m}\left[\frac{1}{\lambda_{\delta}(\infty)}-\frac{1}{\lambda_{\delta}(mN)}\right]\frac{4}{m^2}\sin^2{m\pi\omega} < 0,
\end{align}
which implies that $F^{\infty}<F(N;\delta,\omega)$ for any finite $N\in\mathbb{Z}^+$. Therefore $F(\cdot;\delta,\omega)$ never reach the global minimum at any finite $N$.
\end{proof}

Now we have the following theorem regarding the behavior of $N^*(\gamma)$ as $\gamma\rightarrow\infty$.

\begin{theorem}\label{theorem:CaseIII}
Let $\omega\in(0,1/2)$ and $\delta\in[0,\pi)$ be fixed constants. For the minimization problem (\ref{eqn:min_N}) in which $\lambda_{\delta}(n)$ is the eigenvalue of $\mathcal{L}_{\delta}$ with power kernel (\ref{eqn:Kernel_Power}) and $\alpha\in(\alpha^*, 1)$, the optimal number of bubbles $N^*(\gamma)$, as a function of $\gamma$, is monotonically increasing with no upper bound.
\end{theorem}
\begin{proof}
In this proof we take $N, N^*(\gamma)\in\mathbb{Z}^+$. The case for continuum $N$ is simpler.  The proof of the monotonic increment of $N^*(\gamma)$ follows the same line as that in theorem \ref{theorem:Case1_bound}.

Now we prove by contradiction that there is no upper bound for $N^*(\gamma)$. Assume $\tilde{N}\in\mathbb{Z}^+$ is the optimal upper bound for $N^*(\gamma)$, in other words,
\begin{align*}
    N^*(\gamma) \le \tilde{N},   \text{\ for\ } \forall \gamma \text{\ and\ }  N^*(\gamma) = \tilde{N} \text{\ for\ } \forall \gamma \ge \tilde{\gamma}.
\end{align*}
Note that $F(N;\delta,\omega) \rightarrow F^{\infty}$ and $F^{\infty} < F(N;\delta,\omega)$ for any finite $N\in\mathbb{Z}^+$, it is evident that there exists a $N>\tilde{N}$ such that
\begin{align*}
&E^{\text{tot}}(\tilde{N}, \tilde{\gamma} ) < E^{\text{tot}}(N,\tilde{\gamma}), \\
\text{and}\quad &F(N;\delta,\omega) < F(\tilde{N};\delta,\omega).
\end{align*}
Then it follows that
\begin{align*}
    \tilde{\gamma} \le \frac{\pi(N-\tilde{N})}{F(\tilde{N};\delta,\omega) - F(N;\delta,\omega)} < \bar{\gamma},
\end{align*}
for some $\bar{\gamma}$, which implies
\begin{align}\label{2}
 2\tilde{N} + \frac{2\bar{\gamma}}{\pi} F(\tilde{N};\delta,\omega) >  2N + \frac{2\bar{\gamma}}{\pi} F(N;\delta,\omega),
\end{align}
or equivalently, $E^{\text{tot}}(\tilde{N}, \bar{\gamma} ) > E^{\text{tot}}(N,\bar{\gamma})$. Hence $E^{\text{tot}}$ will not reach the global minimum at $\tilde{N}$ for some $\bar{\gamma}>\tilde{\gamma}$, contradicting to the assumption that $\tilde{N}$ is the optimal upper bound for $N^*(\gamma)$.
\end{proof}

\subsubsection{Nonlocal Effect on the Minimizers for Case IV: $\alpha\in(1,3)$}\label{subsubsection:Case4}
In this subsection, we consider the case in which the eigenvalues $\lambda_{\delta}(n)$ is monotonically increasing as $n\rightarrow\infty$. In this case, it is evident that $F(N;\delta,\omega)$ in (\ref{eqn:F})
\begin{align*}
F(N;\delta,\omega) = \sum_{m=1}^{\infty}\frac{1}{\lambda_{\delta}(mN)}\left(\frac{\sin{m\pi\omega}}{m} \right)^2
\end{align*}
is monotonically decreasing with respect to $N$. Therefore, we have the following theorem, with the same conclusion as theorem \ref{theorem:CaseIII}.
\begin{theorem}\label{theorem:caseIV}
Let $\omega\in(0,1/2)$ and $\delta\in[0,\pi)$ be fixed constants. For the discrete minimization problem (\ref{eqn:min_N}) in which $\lambda_{\delta}(n)$ is the eigenvalue of $\mathcal{L}_{\delta}$ with the power kernel (\ref{eqn:Kernel_Power}) and $\alpha\in(1,3)$, the optimal number of bubbles $N^*(\gamma)$, as a function of $\gamma$, is monotonically increasing with no upper bound.
\end{theorem}

\begin{proof}
First of all, due to the monotonicity of the eigenvalues $\lambda_{\delta}(n)$, it is unnecessary to treat $N$ as a continuum in the proof. We will consider $N, N^*(\gamma) \in \mathbb{Z}^+$.

Secondly, the proof of the monotonic increment of $N^*(\gamma)$ follows the same idea as the proof of monotonic increment in theorem \ref{theorem:CaseIII}.

Lastly, we can prove by contradiction that there is no upper bound for $N^*(\gamma)$, which is similar to that in theorem \ref{theorem:CaseIII}. The only difference is that in the proof of theorem \ref{theorem:CaseIII}, we need to find some $N>\tilde{N}$ such that $F(N;\delta,\omega) < F(\tilde{N};\delta,\omega)$; while for the current proof, due to the monotonic decrement of $F(\cdot;\delta,\omega)$, we have that $F(N;\delta,\omega) < F(\tilde{N};\delta,\omega)$ holds for any $N>\tilde{N}$. We omit the remaining details as they are simply repeat of the proof in theorem \ref{theorem:CaseIII}.
\end{proof}

\begin{remark}
In this case, the optimal number of bubbles $N^*(\gamma)\rightarrow\infty$ as the repulsive strength $\gamma \rightarrow \infty$. This is analogous to the standard OK model in which $\mathcal{L}_{\delta} = -\Delta$.
\end{remark}

\begin{remark}
The above analysis of cases I, III and IV indicates that $N^*(\gamma)$ has a finite upper bound if and only if the eigenvalue $\lambda_{\delta}(n)$ reaches the global maximum at a finite $n$. For the power kernel $\rho_{\delta}(s)$, since $\lambda_{\delta}(n)$ attains the global maximum at a finite $n$ for $\alpha = 0$, then $F(N;\delta,\omega)$ attains the global minimum at a finite $\tilde{N}$, consequently the optimal number of bubbles $N^*$ has an upper bound, $N^*(\gamma) \le \tilde{N}$; on the other hand, when $\alpha\in(\alpha^*,3)$, the supreme of $\lambda_{\delta}(n)$ is attained when $n\rightarrow\infty$, thus the optimal number $N^*$ of bubbles grows unbounded as $\gamma\rightarrow\infty$.
\end{remark}

\subsubsection{$\gamma$-Effect on Minimizers for Case II} 

While it is clear that the boundedness of $N^*(\gamma)$ is determined by the existence of a global maximum of $\lambda_{\delta}(n)$ (and therefore the existence of a global minimum of $F(N;\delta,\omega)$) for cases I, III and IV,   the nonlocal effect on the minimizers for case II is rather subtle.

\begin{figure}[t] 
\begin{center}
\includegraphics[width=0.35\linewidth]{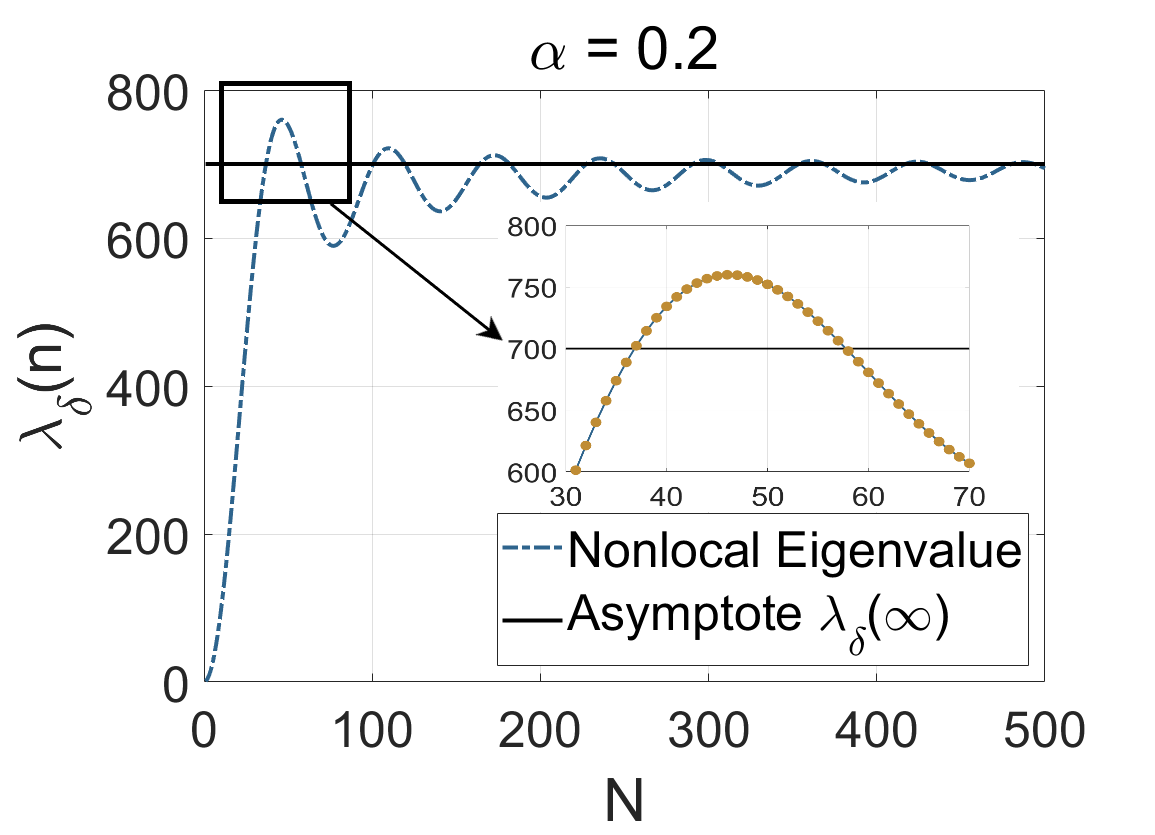}
\includegraphics[width=0.35\linewidth]{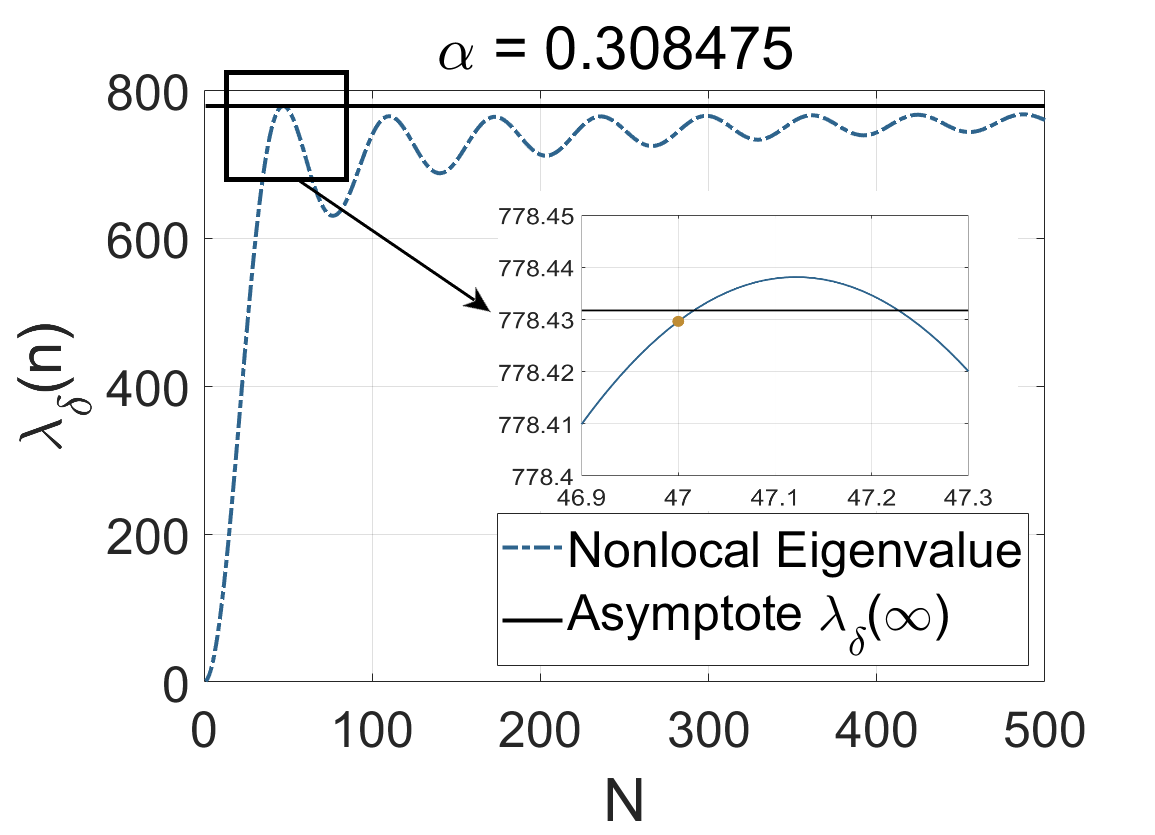}
\end{center}
\caption{Continuous eigenvalue $\lambda_{\delta}(n), n\in\mathbb{R}_+$ (dashed blue curve) v.s. discrete eigenvalue $\lambda_{\delta}(n), n\in\mathbb{Z}_+$ (dark orange dots) for power kernel operator $\rho_{\delta}(s) = \frac{3-\alpha}{\delta^{3-\alpha}|s|^{\alpha}}, s \in [-\delta,0)\cup(0,\delta]$. In left subfigure, both continuous eigenvalue and discrete eigenvalue reach the global maxima somewhere in the first crest.  In the right subfigure, the continuous eigenvalue reaches the global maximum in the first crest, but the discrete eigenvalue does not.}
\label{fig:eigenvalue_powerkernel2} 
\end{figure}

When $\alpha\in(0,\alpha^*)$ and $n\in [1,\infty)$ is treated as a continuum, the eigenvalue $\lambda_{\delta}(n)$ in (\ref{eqn:eigenvalue}) oscillates unevenly around the limit $\lambda_{\delta}(\infty)$ initially, but eventually becomes oscillating from below and converge to $\lambda_{\delta}(\infty)$. More importantly,  $\lambda_{\delta}(n)$ attains the unique global maximum at the first crest. This is indicated by the proof of lemma \ref{lemma:alpha*} and Figure \ref{fig:eigenvalue_powerkernel2}. 

On the other hand, it may lead us to two opposite conclusions when we consider the discrete variable $n\in\mathbb{Z}_+$ in case II. This issue is raised because of the peculiar pattern of $\lambda_{\delta}(n)$. Note that when $\alpha\in(0,\alpha^*)$, there are only a finite number of bounded intervals over which the continuous $\lambda_{\delta}(n)$ is greater than $\lambda_{\delta}(\infty)$. The closer $\alpha$ is to 0, the more number of such bounded intervals and the larger of each interval is. Therefore, when taking $\alpha$ close to 0 and evaluating $\lambda_{\delta}(n)$ over integer grid $n\in\mathbb{Z}_+$, it is highly possible that some integer grid point $\tilde{n}$ falls into one of such bounded intervals so that $\lambda_{\delta}(\tilde{n})>\lambda_{\delta}(\infty)$. Then it follows that a discrete global maximum $n^* = \text{argmax}_{n\in\mathbb{Z}_+} \lambda_{\delta}(n)$ exists, see the left subfigure in Figure \ref{fig:eigenvalue_powerkernel2}. Though it is unclear to us theoretically how $F(N;\delta,\omega), N\in\mathbb{N}_+$ behaves, numerical calculations indicate that $F(N;\delta,\omega)$ reaches a global minimum at some finite $\tilde{N}\in\mathbb{Z}_{+}$, and consequently provides a upper bound for the optimal number $N^{*}(\gamma) \le \tilde{N}$.  On the contrary, when $\alpha$ is sufficiently close to $\alpha^*$, there is only one tiny interval over which the continuous $\lambda_{\delta}(n)$ is greater than $\lambda_{\delta}(\infty)$, see right subfigure in Figure \ref{fig:eigenvalue_powerkernel2}. Such an interval eventually shrinks to a point at the limit of $\alpha\rightarrow\alpha^*$. Therefore when taking $\alpha$ sufficiently close to $\alpha^*$ and evaluating $\lambda_{\delta}(n)$ over integer grid $n\in\mathbb{Z}_+$, it is highly unlikely that any integer grid point will fall into such small interval (but still possible by taking proper values of $\delta$ and $\omega$), leading to the fact that the supremum of $\lambda_{\delta}(n), n\in\mathbb{Z_{+}}$ is attained at $\infty$. Note that this is identical to the case III, in which $N^*(\gamma)\rightarrow\infty$ as $\gamma\rightarrow\infty$.

\subsection{Nonlocal $\gamma$-Effect on Optimizers for Gauss Kernel}\label{subsection:GaussKernel}

In this section, we apply the same argument to the nonlocal effect on the minimizers for Gauss kernel (\ref{eqn:Kernel_Gauss}). Indeed, since the eigenvalue $\lambda_{\delta}(n)$ in (\ref{eqn:Eigenvalue_Gauss}) is monotonically increasing, it matches the condition in case IV for the nonlocal operator with power kernel as discussed in  subsection \ref{subsubsection:Case4}. Consequently, the optimal number of bubbles $N^*(\gamma)$ is unbounded as $\gamma\rightarrow\infty$. We summarize it as below.

\begin{corollary}
Consider the discrete minimization problem (\ref{eqn:min_N}) in which $\lambda_{\delta}(n)$ is given in (\ref{eqn:Eigenvalue_Gauss}) as the eigenvalue of nonlocal operator $\mathcal{L}_{\delta}$ with Gauss-type kernel (\ref{eqn:Kernel_Gauss}). Let $\delta\in(0,\pi)$ and $\omega \in(0, 1/2)$ be fixed constants. The optimal $N^*(\gamma)$, as a function of $\gamma$, is monotonically increasing with no upper bound.
\end{corollary}

\subsection{Nonlocal $\delta$-Effect on Optimizers}\label{subsection:delta}

In this section, we will still consider the minimization problem
\begin{align}\label{eqn:min_main}
\min_{N\in\mathbb{Z}^+} E_0^{\text{NOK}}[U_N^{\omega}] = 2N + \frac{2\gamma}{\pi}\sum_{m=1}^{\infty}\frac{1}{\lambda_{\delta}(mN)}\left(\frac{\sin(m\pi\omega)}{m}\right)^2,
\end{align}
with a focus on the influence of the parameter $\delta$ to the minimizers.  To this end, we denote by $E^{\text{tot}}(N,\gamma,\delta)$ the objective function in (\ref{eqn:min_main}). Our goal is to study the effect of $\delta$ on the optimizer $N^*(\gamma,\delta)$. Note that we set the optimizer $N^*$ as a function of two variables $\gamma$ and $\delta$. In the previous section, we studied the dependence of $N^*$ on $\gamma$ with fixed $\delta$. In this section, on the contrary, we will explore the dependence of $N^*$ on $\delta$ with fixed $\gamma$. For a given nonlocal operator $\mathcal{L}_{\delta}$, if the increment of $\delta$ causes the increment of $N^*$, (treating $N^*$ as a continuum, it means $\frac{\partial N^*}{\partial \delta}\ge0$), we say that it {\it promotes} the bubble splitting, otherwise, it {\it demotes} the bubble splitting.

\subsubsection{Analysis for Monotonically Increasing $\lambda_{\delta}(n)$ }

To ease the analysis, we consider the case in which the eigenvalues $\lambda_{\delta}(n)$ of the nonlocal operator $\mathcal{L}_{\delta}$ are strictly monotonically increasing to $\infty$. In this case, $\lambda_{\delta}^{-1}$ is a strictly decreasing function. We further need $\lambda_{\delta}^{-1}$ to be strictly convex.

\begin{lemma}
Let the eigenvalues $\lambda_{\delta}(n)$ of the nonlocal operator $\mathcal{L}_{\delta}$ be such that $\lambda_{\delta}(n)^{-1}$ is a strictly decreasing and strictly convex function. Then  $F(N;\delta,\omega)$ is a strictly monotonically decreasing and strictly convex function with respect to $N$.
\end{lemma}
\begin{proof}
Given
\[
F(N;\delta,\omega) = \sum_{m=1}^{\infty} \frac{1}{\lambda_{\delta}(mN)} \left(\frac{\sin m\pi\omega}{m}\right)^2,
\]
this is a direct result by verifying the signs of the first and second derivatives of $F$ with respect to $N$.
\end{proof}

Since $F(\cdot; \delta,\omega)$ is strictly decreasing and strictly convex, we have that
\begin{align}\label{eqn:slope_of_energy}
\Delta_N F(N;\delta,\omega): = F(N+1;\delta,\omega)-F(N;\delta,\omega),
\end{align}
is a strictly decreasing function. Therefore the optimal number $N^*$ of bubbles is attained when either one of the following inequalities is uniquely satisfied
\begin{align}
\Delta_N F(N^*-1;\delta,\omega) \le -\frac{\pi}{\gamma} &< \Delta_N F(N^*;\delta,\omega), \quad N^* \ge 2, \\
\text{or}\qquad -\frac{\pi}{\gamma} &< \Delta_N F(N^*;\delta,\omega), \quad N^* = 1.
\end{align}
Consequently, the $\delta$-effect on the promotion/demotion of the bubble splitting is decided by the rate of change of $\Delta_N F$ with respect to $\delta$. More precisely, if 
\[
\frac{\partial}{\partial\delta}\Delta_N F(N;\delta,\omega)\le0 \ (\ge0, \ respectively) 
\]
at $N^*-1$ and $N^*$, the bubble splitting is {\it instantly} promoted (demoted, respectively). If 
\[
\Delta_{\delta}\Delta_N F(N;\delta,\omega): = \Delta_N F(N;\delta,\omega) - \Delta_N F(N;0,\omega) \le 0 \ (\ge 0, \ respectively)
\]
holds at $N^*-1$ and $N^*$, the bubble splitting is {\it cumulatively} promoted (demoted, respectively).

We summarize it as the following theorems.

\begin{theorem}\label{theorem:delta_effect}
Let $\lambda_{\delta}(n)$ be the eigenvalues for the nonlocal operator $\mathcal{L}_{\delta}$ such that $(\lambda_{\delta}(n))^{-1}$ is strictly decreasing and strictly convex. Assuming the following inequalities hold:
\begin{align}\label{eqn:DeltaE_01}
\frac{\partial}{\partial\delta} \Delta_N F(N^*-1;\delta,\omega)\le0, \quad \frac{\partial}{\partial\delta} \Delta_N F(N^*;\delta,\omega) \le 0,
\end{align}
then the nonlocal factor $\delta$ instantly promotes the bubble splitting, namely, $N^*(\gamma,\delta)$ increases with respect to $\delta$. On the contrary, if 
\begin{align}\label{eqn:DeltaE_02}
\frac{\partial}{\partial\delta} \Delta_N F(N^*-1;\delta,\omega)\ge0, \quad \frac{\partial}{\partial\delta} \Delta_N F(N^*;\delta,\omega) \ge 0,
\end{align}
the nonlocal factor $\delta$ instantly demotes the bubble splitting, namely, $N^*(\gamma,\delta)$ decreases with respect to $\delta$.
\end{theorem}

\begin{theorem}\label{theorem:delta_effect2}
Let $\lambda_{\delta}(n)$ be the eigenvalues for the nonlocal operator $\mathcal{L}_{\delta}$ such that $(\lambda_{\delta}(n))^{-1}$ is strictly decreasing and strictly convex. Assuming the following inequalities hold:
\begin{align}\label{eqn:DeltaE_03}
\Delta_{\delta} \Delta_N F(N^*-1;\delta,\omega)\le0, \quad \Delta_{\delta} \Delta_N F(N^*;\delta,\omega) \le 0,
\end{align}
then the nonlocal factor $\delta$ cumulatively promotes the bubble splitting, namely, $N^*(\gamma,\delta)\ge N^*(\gamma,0)$. On the contrary, if 
\begin{align}\label{eqn:DeltaE_04}
\Delta_{\delta}  \Delta_N F(N^*-1;\delta,\omega)\ge0, \quad \Delta_{\delta}  \Delta_N F(N^*;\delta,\omega) \ge 0,
\end{align}
the nonlocal factor $\delta$ cumulatively demotes the bubble splitting, namely, $N^*(\gamma,\delta) \le N^*(\gamma,0)$.
\end{theorem}

\begin{remark}
In regard to theorem \ref{theorem:delta_effect}, it is evident that if
\begin{align}
\frac{\partial}{\partial\delta} \Delta_N F(N^*-1;\delta,\omega)\le0, \quad \frac{\partial}{\partial\delta} \Delta_N F(N^*;\delta,\omega) \ge 0,
\end{align}
then $N^*(\gamma,\delta)$ is unchanged. However, if 
\begin{align}
\frac{\partial}{\partial\delta} \Delta_N F(N^*-1;\delta,\omega)\ge0, \quad \frac{\partial}{\partial\delta} \Delta_N F(N^*;\delta,\omega) \le 0,
\end{align}
it is undetermined for the change of $N^*$. However, it is very rare for this undetermined case to occur as we can take sufficiently small or large values of $\gamma$ such that $\frac{\partial}{\partial\delta} \Delta_N F$ has the same sign at $N^*-1$ and $N^*$. It is similar for the cumulative case in theorem \ref{theorem:delta_effect2}. See the case studies in subsection \ref{subsubsection:powerkernel} and \ref{subsubsection:fractionalLaplacian} for the details.
\end{remark}

\begin{remark}
In this section, we only consider the case in which $\lambda_{\delta}^{-1}$ is strictly decreasing and strictly convex so that the optimizer $N^*(\gamma,\delta)$ is uniquely determined. For other more complicated cases, we can still verify the nonlocal effect on the promotion/demotion of the bubble splitting by numerical methods. However it is lack of analysis tools to perform theoretical study on this subject.
\end{remark}

Though conditions (\ref{eqn:DeltaE_01}) and (\ref{eqn:DeltaE_02}) for instant promotion/demotion and  (\ref{eqn:DeltaE_03}) and (\ref{eqn:DeltaE_04}) for cumulative promotion/demotion are only sufficient, they are still of practical importance as demostrated by the three case studies in the following subsections.

\subsubsection{A Case Study: Power Kernel with $\alpha\in(2,3)$} \label{subsubsection:powerkernel}

Take the nonlocal operator $\mathcal{L}_{\delta}$ with power kernel (\ref{eqn:Kernel_Power}) as an example. The eigenvalues $\lambda_{\delta}(n)$ are given as in (\ref{eqn:eigenvalue}). To apply theorem \ref{theorem:delta_effect} to this example, we first show that $\lambda_{\delta}^{-1}$ is strictly decreasing and strictly convex, 
then explore the condition under which (\ref{eqn:DeltaE_01}) holds.

\begin{lemma}
The eigenvalues $\lambda_{\delta}(n)$ for the nonlocal operator $\mathcal{L}_{\delta}$ with power kernel (\ref{eqn:Kernel_Power}) and $\alpha\in(2,3)$ satisfies
\[
\frac{\partial}{\partial n} \lambda_{\delta}^{-1} < 0, \quad \frac{\partial^2}{\partial n^2} \lambda_{\delta}^{-1} > 0,
\]
namely, $\lambda_{\delta}^{-1}$ is strictly decreasing and strictly convex.
\end{lemma}

\begin{proof}
The eigenvalue $\lambda_{\delta}$ for the nonlocal operator $\mathcal{L}_{\delta}$ with power kernel $\rho(\xi) = \frac{3-\alpha}{\xi^{\alpha}}, \xi\in(0,1)$ and $\alpha\in(2,3)$ reads
\[
\lambda_{\delta} = \frac{2(3-\alpha)}{\delta^2} \int_0^1 \frac{1-\cos n\delta\xi}{\xi^{\alpha}} d\xi : = \frac{2(3-\alpha)}{\delta^2} K(n\delta).
\]
Taking the first order derivative, we have
\[
\dfrac{\partial \lambda_{\delta}}{\partial n} = \frac{2(3-\alpha)}{\delta} \int_0^1 \frac{\sin n\delta\xi}{\xi^{\alpha-1}} d\xi > 0, \quad \forall n\delta > 0.
\]
Hence $\lambda_{\delta}^{-1}(n)$ is strictly decreasing.

Taking the second order derivative for $\lambda_{\delta}^{-1}$, we get
\[
\frac{\partial^2}{\partial n^2} \lambda_{\delta}^{-1} = \lambda_{\delta}^{-3}\bigg[2\bigg(\frac{\partial \lambda_{\delta}}{\partial n}\bigg)^2 - \lambda_{\delta}\frac{\partial^2\lambda_{\delta}}{\partial n^2}\bigg]
=  \lambda_{\delta}^{-3}\frac{4(3-\alpha)^2}{\delta^2}\Big[2\bigg(K'(n\delta)\Big)^2 - K(n\delta)K''(n\delta)\bigg].
\]
It suffices to show that 
\begin{align}\label{eqn:K_ineqn}
\Big[2\bigg(K'(n\delta)\Big)^2 - K(n\delta)K''(n\delta)\bigg] > 0.
\end{align}
Note that $K(\cdot)$ satisfies the following differential equation \cite{DuYang_JCP2017}
\begin{align}\label{eqn:Kp}
K'(x) = \frac{\alpha - 1}{x} K(x) + \frac{1-\cos x}{x}.
\end{align}
Differentiating the above identity (\ref{eqn:Kp}) yields,
\begin{align}\label{eqn:Kpp}
K''(x) = \frac{\alpha - 2}{x} K'(x) + \frac{\sin x}{x}.
\end{align}
Using the equation (\ref{eqn:Kp}) to replace $K$ by $K'$ and the equation (\ref{eqn:Kpp}) to replace $K''$ by $K'$, the left hand side of the equation (\ref{eqn:K_ineqn}), after simplification and multiplication by $(\alpha-1)$, becomes
\begin{align*}
&\alpha \Big( K'(n\delta) \Big)^2 - \bigg(\sin n\delta - (\alpha-2)\frac{1-\cos n\delta}{n\delta}\bigg) K'(n\delta) + \sin n\delta \cdot \frac{1-\cos n\delta}{n\delta} \\
> & \Big( K'(n\delta) \Big)^2 - \bigg(\sin n\delta + \frac{1-\cos n\delta}{n\delta}\bigg) K'(n\delta) + \sin n\delta \cdot \frac{1-\cos n\delta}{n\delta} \\
= & \bigg(K'(n\delta) - \sin n\delta\bigg)\bigg( K'(n\delta) - \frac{1-\cos n\delta}{n\delta} \bigg) = \text{I}\cdot \text{II}, \quad \forall \alpha\in(2,3).
\end{align*}

It remains to show that $\text{I}, \text{II} \ge 0$. When $n\delta \le \frac{2\pi}{3}$, $\frac{\sin n\delta\xi}{\xi^{\alpha-1}}$ is decreasing with respect to $\xi$ over (0,1], then
\begin{align}\label{eqn:Kp_lowerbound1}
K'(n\delta) = \int_0^1 \frac{\sin n\delta\xi}{\xi^{\alpha-1}} d\xi \ge \int_0^1 \frac{\sin n\delta}{1^{\alpha-1}} d\xi = \sin n\delta.
\end{align}
Besides, since $\sin x \ge \frac{1-\cos x}{x}$ for $x\in[0,\frac{2\pi}{3}]$, we also have
\begin{align}\label{eqn:Kp_lowerbound2}
K'(n\delta) \ge \sin n\delta \ge \frac{1 - \cos n\delta}{n\delta} .
\end{align}
When $n\delta > \frac{2\pi}{3}$ and $\alpha\in(2,3)$,
\begin{align}\label{eqn:Kp_lowerbound}
K'(n\delta) = (n\delta)^{\alpha-2} \int_0^{n\delta} \frac{\sin t}{t^{\alpha-1}} dt > \left(\frac{2\pi}{3}\right)^{\alpha-2}\int_0^{n\delta} \frac{\sin t}{t} dt > \left(1\right)^{\alpha-2}\int_0^{2\pi} \frac{\sin t}{t} dt \approx 1.42,
\end{align}
in which the first inequality is due to the fact that $\int_0^{n\delta} \frac{\sin t}{t^{\alpha-1}} dt$ is an increasing function with respect to $\alpha \in (2,3)$, and the second inequality is due to the fact that the sine integral $\text{Si}(x) = \int_0^x \frac{\sin t}{t} dt$ reaches the unique global minimum at $x = 2\pi$ over the interval $x\in(\frac{2\pi}{3}, \infty)$.
Then it follows that
\begin{align}\label{eqn:Kp_lowerbound3}
K'(n\delta) \ge 1 \ge \sin n\delta.
\end{align}
Additionally, note that $\frac{1 - \cos n\delta}{n\delta} < \frac{3}{\pi}$ for $n\delta > \frac{2\pi}{3}$, we can continue with the estimate (\ref{eqn:Kp_lowerbound}) to get
\begin{align}\label{eqn:Kp_lowerbound4}
K'(n\delta) > \left(\frac{2\pi} {3}\right)^{\alpha-2}\int_0^{2\pi} \frac{\sin t}{t} dt > \frac{3}{\pi} > \frac{1 - \cos n\delta}{n\delta}.
\end{align}
Combining (\ref{eqn:Kp_lowerbound1}), (\ref{eqn:Kp_lowerbound2}), (\ref{eqn:Kp_lowerbound3}) and (\ref{eqn:Kp_lowerbound4}), we get $\text{I}, \text{II} \ge 0$. Consequently (\ref{eqn:K_ineqn}) holds, and therefore $\frac{\partial^2}{\partial n^2} \lambda_{\delta}^{-1} > 0$. 
\end{proof}

Now we study the condition under which (\ref{eqn:DeltaE_01}) holds. Note that
\begin{align}\label{eqn:series_delta}
\frac{\partial}{\partial \delta}\Delta_N F(N;\delta,\omega) = 
\sum_{m=1}^{\infty}\frac{\partial}{\partial \delta}\left(\frac{1}{\lambda_{\delta}(m(N+1))}-\frac{1}{\lambda_{\delta}(mN)}\right)\left(\frac{\sin{m\pi\omega}}{m}\right)^2,
\end{align}
the sign of $\frac{\partial}{\partial \delta}\Delta_N F(N;\delta,\omega)$ is determined by the sign of 
\[
\frac{\partial}{\partial \delta}\left(\frac{1}{\lambda_{\delta}(m(N+1))}-\frac{1}{\lambda_{\delta}(mN)}\right), \quad m=1,2,\cdots.
\]
Differentiating $\lambda_{\delta}^{-1}(n)$ with respect to $\delta$, we have
\begin{align*}
\frac{\partial}{\partial \delta} \frac{1}{\lambda_{\delta}(n)} &= \frac{\partial}{\partial \delta} \left( \frac{\delta^2}{2(3-\alpha)}\frac{1}{K(n\delta)}  \right) \\
&= \frac{\delta}{3-\alpha}\frac{1}{K(n\delta)} - \frac{\delta}{2(3-\alpha)}\frac{1}{K(n\delta)^2}K'(n\delta)(n\delta) \\
& = \frac{\delta}{3-\alpha}\frac{1}{K(n\delta)} - \frac{\delta}{2(3-\alpha)}\frac{1}{K(n\delta)^2}\Big((\alpha-1)K(n\delta)+(1-\cos n\delta)\Big) \\
& = \frac{3-\alpha}{\delta} \frac{1}{\lambda_{\delta}(n)} - \frac{2(3-\alpha)}{\delta^3} \frac{1-\cos n\delta}{\lambda_{\delta}(n)^2}.
\end{align*}
then
\begin{align}
&\frac{\partial}{\partial \delta}\left(\frac{1}{\lambda_{\delta}(m(N+1))}-\frac{1}{\lambda_{\delta}(mN)}\right) \nonumber\\
=& \frac{3-\alpha}{\delta} \left(\frac{1}{\lambda_{\delta}(m(N+1))}-\frac{1}{\lambda_{\delta}(mN)}\right)
- \frac{2(3-\alpha)}{\delta^3} \left(\frac{1-\cos m(N+1)\delta}{\lambda_{\delta}(m(N+1))^2}-\frac{1-\cos mN\delta}{\lambda_{\delta}(mN)^2}\right). \label{eqn:temp_delta}
\end{align}
Note that \cite{DuYang_SINA2016} for $\alpha\in(2,3)$,
\begin{align}\label{eqn:bound_lambda}
C_1(\delta) n^{\alpha-1}  \le \lambda_{\delta}(n) \le C_2(\delta) n^{\alpha-1},
\end{align}
which, together with equation (\ref{eqn:temp_delta}), implies that $\frac{\partial}{\partial \delta}\left(\frac{1}{\lambda_{\delta}(m(N+1))}-\frac{1}{\lambda_{\delta}(mN)}\right)$ is uniformly bounded with respect to $m$. Therefore the series in (\ref{eqn:series_delta}) is absolutely convergent. 

Now we show that for $\alpha \in (2,3)$, $D_1(\delta) n^{\alpha-2} \le \frac{\partial}{\partial n}\lambda_{\delta}(n) \le D_2(\delta) n^{\alpha-2}$. Note that
\[
\dfrac{\partial \lambda_{\delta}}{\partial n} = \frac{2(3-\alpha)}{\delta^{3-\alpha}} n^{\alpha-2} \int_0^{n\delta} \frac{\sin s}{s^{\alpha-1}} ds : =  \frac{2(3-\alpha)}{\delta^{3-\alpha}} n^{\alpha-2} Q_n,
\]
it remains to prove that the positive sequence $\{Q_n\}$ has uniform upper and lower bounds independent of $n$. For the upper bound, 
\[
Q_n = \int_0^{n\delta} \frac{\sin s}{s^{\alpha-1}} ds \le \int_0^{\pi} \frac{\sin s}{s^{\alpha-1}} ds.
\]
For the lower bound, we know that if $n\delta \ge \pi$, 
\[
Q_n = \int_0^{n\delta} \frac{\sin s}{s^{\alpha-1}} ds \ge \int_0^{2\pi} \frac{\sin s}{s^{\alpha-1}} ds;
\]
and if $n\delta < \pi$, 
\[
Q_n = \int_0^{n\delta} \frac{\sin s}{s^{\alpha-1}} ds \ge \int_0^{\delta} \frac{\sin s}{s^{\alpha-1}} ds,
\]
therefore
\[
Q_n  \ge \min \left\{  \int_0^{2\pi} \frac{\sin s}{s^{\alpha-1}} ds, \int_0^{\delta} \frac{\sin s}{s^{\alpha-1}} ds \right\},
\]
which indicates the bound
\begin{align}\label{eqn:bound_lambda_prime}
D_1(\delta) n^{\alpha-2} \le \frac{\partial}{\partial n}\lambda_{\delta}(n) \le D_2(\delta) n^{\alpha-2}.
\end{align}

Note that $\frac{\partial}{\partial n} \lambda_{\delta}^{-1}(n) = -\lambda_{\delta}^{-2}(n) \frac{\partial}{\partial n}\lambda_{\delta}(n)$, using the bounds (\ref{eqn:bound_lambda}) and (\ref{eqn:bound_lambda_prime}), we have
\begin{align}
-C_1^{-2}D_1 n^{-\alpha} \le \frac{\partial}{\partial n} \lambda_{\delta}^{-1}(n) \le -C_2^{-2}D_2 n^{-\alpha}.
\end{align}
Taking integral over $[mN, m(N+1)]$, we can easily obtain the following bounds,
\begin{align}
-\frac{C_1^{-2}D_1}{(\alpha-1)N^{\alpha-1}}\frac{1}{m^{\alpha-1}} \le \frac{1}{\lambda_{\delta}(m(N+1))}-\frac{1}{\lambda_{\delta}(mN)} \le -\frac{C_2^{-2}D_2}{(\alpha-1)(N+1)^{\alpha}}\frac{1}{m^{\alpha-1}}.
\end{align}
By (\ref{eqn:temp_delta}), we have
\begin{align*}
\frac{\partial}{\partial \delta}\left(\frac{1}{\lambda_{\delta}(m(N+1))}-\frac{1}{\lambda_{\delta}(mN)}\right) & \le
\frac{3-\alpha}{\delta} \left(\frac{1}{\lambda_{\delta}(m(N+1))}-\frac{1}{\lambda_{\delta}(mN)}\right) + \frac{2(3-\alpha)}{\delta^3} \frac{2}{\lambda_{\delta}(mN)^2} \\
& \le - \frac{3-\alpha}{\delta} \frac{C_2^{-1}D_2}{(\alpha-1)(N+1)^{\alpha}}\frac{1}{m^{\alpha-1}} + 
 \frac{4(3-\alpha)}{\delta^3} \frac{C_1^{-2}}{N^{2\alpha-2}}\frac{1}{m^{2\alpha-2}} ,
\end{align*}
Note that the sign of the right hand size of the above inequality is controlled by $\frac{1}{(N+1)^{\alpha}}$ and $\frac{1}{N^{2\alpha-2}}$. Since $\frac{1}{(N+1)^{\alpha}}$ decays slower than $\frac{1}{N^{2\alpha-2}}$ as $N\rightarrow\infty$, it follows that for a sufficiently large but fixed $N$,
\begin{align}\label{inequality:lambda_delta}
\frac{\partial}{\partial \delta}\left(\frac{1}{\lambda_{\delta}(m(N+1))}-\frac{1}{\lambda_{\delta}(mN)}\right) \le 0, \quad \text{for\ } m = 1,2,\cdots.
\end{align}

Therefore, we have the following theorem hold,

\begin{theorem}
Let the nonlocal operator $\mathcal{L}_{\delta}$ be with the power kernel (\ref{eqn:Kernel_Power}) and $\alpha\in(2,3)$, and assme the optimizer $N^* = N^*(\gamma,\delta)$ is sufficiently large, then the inequalities in (\ref{eqn:DeltaE_01}) hold, and the nonlocal factor $\delta$ instantly promotes the bubble splitting.
\end{theorem}

\begin{proof}
By theorem \ref{theorem:caseIV}, we know that for a fixed $\delta$, 
\[
N^*(\gamma,\delta)\rightarrow\infty,\quad \text{as\ } \gamma\rightarrow\infty.
\]
Hence we can take sufficiently large $\gamma$ such that $N^*(\gamma,\delta)$ is sufficiently large and the inequality (\ref{inequality:lambda_delta}) holds for any $m\in\mathbb{Z}_+$, leading to $\frac{\partial}{\partial \delta}\Delta_N F(N;\delta,\omega)\le 0$ by (\ref{eqn:series_delta}). Therefore, the conditions in theorem \ref{theorem:delta_effect} hold, implying the instant promotion of the bubble splitting by $\delta$.
\end{proof}

\subsubsection{A Case Study: Gauss-type Kernel} \label{subsubsection:fractionalLaplacian}

The analysis of $\delta$-effect for the nonlocal operator $\mathcal{L}_{\delta}$ with Gauss-type kernel is rather straightforward in comparison with that for the power kernel.

Note that the nonlocal operator $\mathcal{L}_{\delta}$ with Gauss-type kernel (\ref{eqn:Kernel_Gauss}) has eigenvalues $\lambda_{\delta}(n)$ in simple closed form (\ref{eqn:Eigenvalue_Gauss}), it is evident to verify that $\lambda_{\delta}^{-1}$ is strictly decreasing and strictly convex.

\begin{lemma}
The eigenvalues $\lambda_{\delta}(n)$ of the nonlocal operator $\mathcal{L}_{\delta}$ with Gauss-type kernel (\ref{eqn:Kernel_Gauss}) satisfy
\begin{align*}
    \frac{\partial}{\partial n}\lambda_{\delta}^{-1}(n) < 0, \ \ \frac{\partial^2}{\partial n^2}\lambda_{\delta}^{-1}(n)>0.
\end{align*}
\end{lemma}

Now we verify the conditions in theorems \ref{theorem:delta_effect} and \ref{theorem:delta_effect2}. Similar to the case study for the power kernel, the sign of $\frac{\partial}{\partial \delta}\Delta_{N}F(N;s,\omega)$ (and $\Delta_{\delta}\Delta_{N}F(N;s,\omega)$, respectively), is determined by the sign of  
\begin{align*}
   & \frac{\partial}{\partial \delta}\left(\frac{1}{\lambda_{\delta}(m(N+1))}-\frac{1}{\lambda_{\delta}(mN)}\right), \ \ \forall m\in\mathbb{Z}_+, \\
& \Bigg( \text{and}\  \left(\frac{1}{\lambda_{\delta}(m(N+1))}-\frac{1}{\lambda_{\delta}(mN)}\right) -  \left(\frac{1}{\lambda_{0}(m(N+1))}-\frac{1}{\lambda_{0}(mN)}\right), \ \ \forall m\in\mathbb{Z}_+,\ \text{respectively} \Bigg).
\end{align*}
In this case, we claim that $\frac{\partial}{\partial n}\left(\frac{\partial}{\partial \delta}\frac{1}{\lambda_{\delta}(n)}\right) \ge 0$. Indeed, by tedious calculation, we have that
\begin{align*}
\frac{\partial}{\partial n}\left(\frac{\partial}{\partial \delta}\lambda_{\delta}^{-1}(n)\right) =  \dfrac{n\delta^3 e^{v}}{4\left(e^{v}-1\right)^3}\left[e^v(v-2)+ v + 2 \right]\geq 0
\end{align*}
in which $v = (n\delta)^2/4$. Consequently both inequalities (\ref{eqn:DeltaE_02}) and (\ref{eqn:DeltaE_04}) hold, which implies, by theorems \ref{theorem:delta_effect} and \ref{theorem:delta_effect2}, both instant demotion and cumulative demotion of the bubble splitting induced by $\delta$. 

We summarize the conclusion as the following corollary.
\begin{corollary}
Let the nonlocal operator $\mathcal{L}_{\delta}$ be with the Gauss-type kernel (\ref{eqn:Kernel_Gauss}), then the nonlocal factor $\delta$ both instantly and cumulatively demotes the bubble splitting. 
\end{corollary}

\subsubsection{A Case Study: $(\delta I - \Delta)^{-1}$}

In the two case studies above, we note that it is the behaviors of $\partial_n \lambda_{\delta}^{-1}, \partial_{n}\partial_{n} \lambda_{\delta}^{-1}$ and $\partial_{n}\partial_{\delta} \lambda_{\delta}^{-1}$ that determine the $\delta$-effect on demotion/promotion of the bubble splitting. Thereby the analysis framework is applicable not only for the nonlocal operator $\mathcal{\delta}$ but also for some local operators involving some system parameter $\delta$. One typical example is the screened Poisson operator $\mathcal{L}_{\delta} = \delta I - \Delta$ \cite{RenTruskinovsky_Elasticity2000, XuZhao_JSC2020}. Note that here $\delta$ does not represent the nonlocal horizon parameter, but rather a ``screening" constant such that $\mathcal{L}_{\delta} \rightarrow -\Delta$ as $\delta\rightarrow 0$.

The screened Poisson operator 
\begin{align}\label{eqn:ScreenPoisson}
\mathcal{L}_{\delta} = \delta I - \Delta
\end{align}
under periodic boundary condition is associated with eigenvalues
\begin{align}\label{eig:ScreenPoisson}
\lambda_{\delta}(n) = \delta + n^2.
\end{align}
By simple calculation, it follows that
\begin{align}
&\frac{\partial}{\partial n} \lambda_{\delta}^{-1} = \frac{-2n}{(\delta+n^2)^2}<0, \\
&\frac{\partial^2}{\partial n^2} \lambda_{\delta}^{-1} = \frac{2(3n^2-\delta)}{(\delta+n^2)^3} > 0, \quad \text{for\ } n\ge 1, \delta < 3, \\
&\frac{\partial}{\partial n} \frac{\partial}{\partial \delta} \lambda_{\delta}^{-1} = \frac{4n}{(\delta+n^2)^3} > 0.
\end{align}
Consequently the screening constant $\delta$ demotes the bubble splitting, instantly and cumulatively. 
\begin{corollary}
Consider the minimization problem (\ref{eqn:min_N}) with eigenvalues given in (\ref{eig:ScreenPoisson}) for the screened Poisson operator (\ref{eqn:ScreenPoisson}). If $\delta\in(0,3)$, the screening parameter $\delta$ both instantly and cumulatively demotes the bubble splitting.
\end{corollary}

\begin{figure}[t] 
\begin{center}
\includegraphics[width=0.5\linewidth]{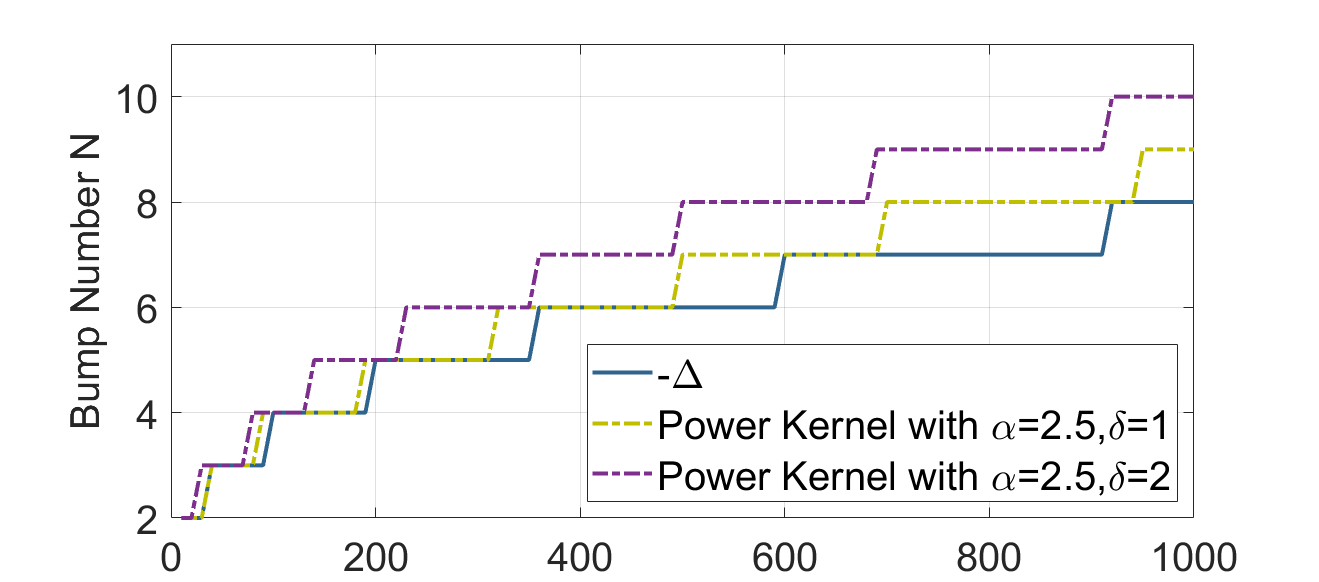}
\includegraphics[width=0.5\linewidth]{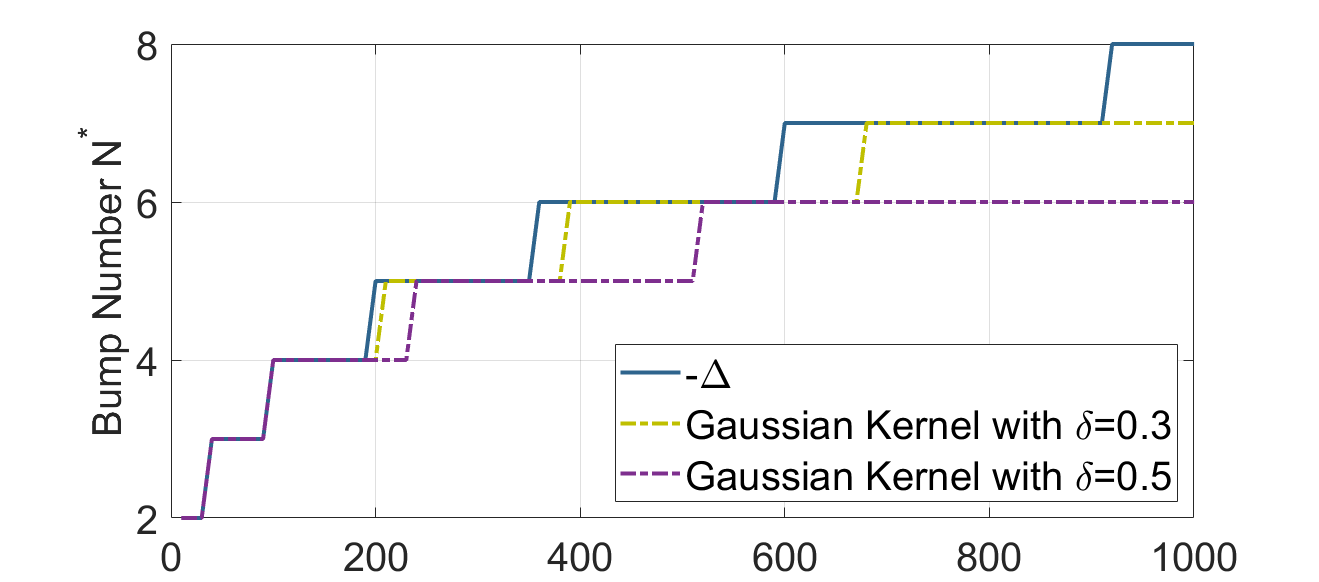}
\includegraphics[width=0.5\linewidth]{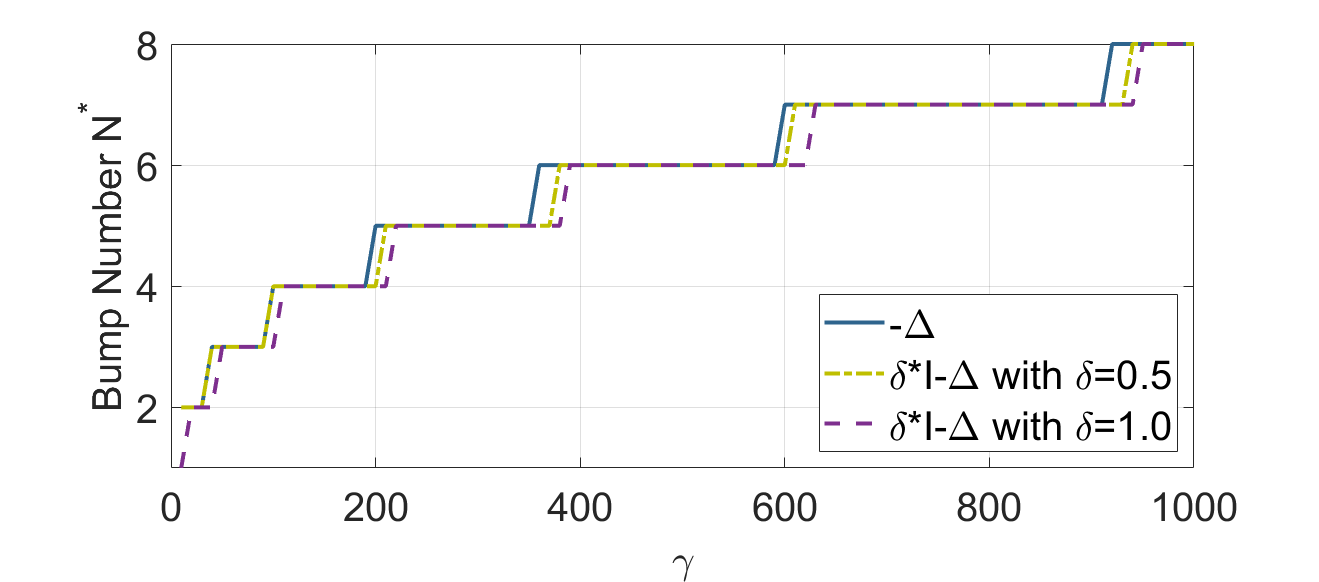}
\end{center}
\caption{Numerical validation for $\delta$-effect on promotion/demotion of bubble splitting. Top: $\delta$ promotes the bubble splitting for power kernel with $\alpha = 2.5$ fixed. Middle: $\delta$ demotes the bubble splitting for Gaussian kernel. Bottom: $\delta$ demotes the bubble splitting for screened Poisson operator. The relative volume $\omega = 0.3$ is fixed.}
\label{fig:promotion} 
\end{figure}

In Figure \ref{fig:promotion}, we numerically verify the $\delta$-effect on promotion/demotion of bubble splitting. In the numerical simulations, the eigenvalues (\ref{eqn:eigenvalue}) corresponding to power kernel are evaluated by a built-in function \texttt{integral} in Matlab. The eigenvalues for Gauss kernel and screened Poisson operator are analytically available. The long-range interaction $F(N;\delta,\omega)$ is approximated by taking the truncated Fourier series with $N_{\text{sum}} = 500$ terms. The relative volume $\omega = 0.3$ is fixed. In the top subfigure, it shows that $\delta$ can promote the bubble splitting. Though we only prove that $\delta$ instantly promotes the bubble splitting in subsection \ref{subsubsection:powerkernel}, the numerical result supports the cumulative promotion as well. The middle (respectively, bottom) subfigure presents the demotion effect of $\delta$ on bubble splitting for Gauss kernel (respectively, screened Poisson operator). Interestingly, when comparing the middle and bottom subfigures, we notice that $\delta$-effect on bubble demotion is stronger for Gauss kernel than that for screened Poisson operator in the sense that for sufficiently large $\gamma$, there exist wider intervals in which $N^*_{\text{Gauss}} < N^*_{-\Delta}$; on the other hand, the intervals of $\gamma$ in which $N^*_{\text{screened}} < N^*_{-\Delta}$ are rather narrower.

\section{Concluding Remark}\label{section:Conclusion}

In this paper, we consider a generalized Ohta-Kawasaki model in which the long-range interaction is induced by a nonlocal operator $\mathcal{L}_{\delta}$. Under periodic boundary condition and the assumption (\ref{eqn:conjucture}) that $U_N^{\omega}$, the $N$-bubble function of equal size and equal distance, minimizes $E_0^{\text{NOK}}$ over $\mathcal{A}_N^{\omega}$,  we study the minimization problem (\ref{eqn:min_N}), minimizing $E_0^{\text{NOK}}[U_N^{\omega}]$ over $N\in\mathbb{Z}_+$. We find that the optimal number $N^*(\gamma)$ of (\ref{eqn:min_N}), may or may not have an upper bound, depending on whether $F(N;\delta,\omega)$ in (\ref{eqn:F}) has a lower bound or not. Additionally, we explore the dependence of $N^*$ on $\delta$. Under some mild conditions on $\lambda_{\delta}$ and $F(N; \delta, \omega)$, we show that nonlocal operator with power kernel ($\alpha\in(2,3)$) promotes the bubble splitting instantly; while the nonlocal operator with Gauss kernel and screened Poisson operator demote the bubble splitting instantly and cumulatively. 

While it is still unclear whether it can be theoretically proved that $U_{N}^{\omega}$ minimizes $E_0^{\text{NOK}}$ over $\mathcal{A}_N^{\omega}$, which will be studied in the future work, several efforts are being made for the ongoing work. Firstly, we design some asymptotically compatible and energy stable numerical schemes to solve the $L^2$ gradient flow dynamics for NOK model. When applying to the two-dimensional case, the NOK system displays some new patterns such as square lattice pattern and hexagonal pattern with elliptical bubbles. Secondly, we are also interested in the nonlocal effect on generalized (ternary) Ohta-Nakazawa model. This study will be even more complicated as there are four parameters $[\gamma_{ij}]_{2\times 2}$ for the long-range interaction, and the patterns of the minimizers of $E_0^{\text{NON}}$ will be more diverse. Thirdly, in the current work, we only introduce a nonlocal operator $\mathcal{L}_{\delta}$ in the long-range interaction term. In the future work, we will consider a more general free energy functional
\begin{align*}
E^{\text{NOK}}_{\epsilon}[u] = \int_{\mathbb{T}^d} \dfrac{\epsilon}{2}|\mathcal{L}_{\delta}^{\frac{1}{2}} u|^2 + \dfrac{1}{\epsilon}W(u)\ \text{d}x + \dfrac{\gamma}{2}\int_{\mathbb{T}^d} |(\mathcal{L}_{\delta})^{-\frac{1}{2}}(f(u)-\omega)|^2\ \text{d}x,
\end{align*}
and study its $\Gamma$-limit, the structure of minimizers, and the nonlocal effect of $\gamma$ and $\delta$ on the minimizers.

\section{Acknowledgements}

Y. Zhao's work is supported by a grant from the Simons Foundation through Grant No. 357963 and NSF grant DMS-2142500.


\newpage


\bibliography{OhtaKawasaki}

\end{document}